\documentclass[11pt]{article}
\usepackage{amssymb, amsthm, amsmath, amscd}
\newtheorem{thm}{Theorem}[section]
\newtheorem{lem}[thm]{Lemma}
\newtheorem{prop}[thm]{Proposition}
\newtheorem{cor}[thm]{Corollary}
\theoremstyle{definition}\newtheorem{df}[thm]{Definition}
\theoremstyle{definition}\newtheorem{rem}[thm]{Remark}
\theoremstyle{definition}\newtheorem{exm}[thm]{Example}

\renewcommand{\phi}{\varphi}

\newcommand{\N}{\mathbb{N}}
\newcommand{\Z}{\mathbb{Z}}
\newcommand{\Q}{\mathbb{Q}}
\newcommand{\R}{\mathbb{R}}
\newcommand{\C}{\mathbb{C}}
\newcommand{\T}{\mathbb{T}}

\newcommand{\Homeo}{\operatorname{Homeo}}
\newcommand{\Aut}{\operatorname{Aut}}

\newcommand{\Aff}{\operatorname{Aff}}

\newcommand{\Isom}{\operatorname{Isom}}
\newcommand{\Tr}{\operatorname{Tr}}
\newcommand{\id}{\operatorname{id}}
\newcommand{\Ima}{\operatorname{Im}}
\newcommand{\Coker}{\operatorname{Coker}}
\newcommand{\Ker}{\operatorname{Ker}}

\newcommand{\xa}{(X,\alpha)}
\newcommand{\yb}{(Y,\beta)}

\newcommand{\CA}{$C^*$-algebra}

\title{Minimal dynamical systems on the product \\
of the Cantor set and the circle}
\author{Huaxin Lin and Hiroki Matui}
\date{}

\begin{document}
\maketitle

\begin{abstract}
We prove that a crossed product algebra arising from a minimal
dynamical system on the product of the Cantor set and the circle
has real rank zero if and only if that system is rigid. In the
case that cocycles take values in the rotation group, it is also
shown that rigidity implies tracial rank zero, and in particular,
the crossed product algebra is isomorphic to a unital simple
AT-algebra of real rank zero. Under the same assumption, we show
that two systems are approximately $K$-conjugate if and only if
there exists a sequence of isomorphisms between two associated
crossed products which approximately maps $C(X\times \T)$
onto $C(X\times \T)$.
\end{abstract}

\section{Introduction}

A celebrated theorem of Giordano, Putnam and Skau \cite{GPS} gave
a dynamical characterization of isomorphism of the crossed product
$C^*$-algebras arising from minimal dynamical systems on the
Cantor set. The $C^*$-algebra theoretic aspect of this result is
indebted to the fact that the algebras are unital simple AT-algebras
with real rank zero. From the work of Q. Lin and N. C. Phillips
\cite{LP2} and the classification of simple \CA s of
tracial rank zero (see \cite{Lnan} and \cite{L3}), crossed product
\CA s arising from minimal diffeomorphisms on a manifold are
isomorphic if they have the same Elliott invariants in the case
that they have real rank zero (see also \cite{LP1}). However,
there are no longer any dynamical characterizations of isomorphism
of these algebras. The notion of approximate conjugacy was
introduced in \cite{LM} where it was also suggested that certain
approximate version of conjugacy may be a right equivalence
relation to ensure isomorphism of crossed product $C^*$-algebras.
Indeed, a complete description was given for Cantor minimal
systems. In the present paper, we consider minimal dynamical
systems on the product of the Cantor set $X$ and the circle $\T$,
and analyze the associated crossed product $C^*$-algebras and
approximate conjugacy for those systems.

Since the Cantor set $X$ is totally disconnected and the circle
$\T$ is connected, every minimal dynamical system on $X\times\T$
can be viewed as a skew product extension of a minimal dynamical
system on $X$. This observation enables us to find a nice
``large'' subalgebra of the crossed product $C^*$-algebra in a
similar fashion to the Cantor case. By applying results of
\cite{Ph2}, it will be shown that the algebra has stable rank one
and satisfies Blackadar's fundamental comparability property.
Moreover, we also prove that the algebra has real rank zero if and
only if every invariant measure on $X$ uniquely extends to an
invariant measure on $X\times\T$.
Such a system is said to be rigid.
As a special case we consider cocycles taking their values
in the rotation group and show that the algebra has tracial rank zero
if and only if the system is rigid.
By an easy computation of $K$-groups and the classification
theorem of \cite{L3}, we conclude that the associated crossed
product $C^*$-algebras are actually unital simple AT-algebras of real
rank zero.

A natural definition of approximate conjugacy for two dynamical
systems $\xa$ and $\yb$ is the following: there exists a sequence
of homeomorphisms $\sigma_n:X\rightarrow Y$ such that
$\sigma_n\alpha\sigma_n^{-1}$ converges to $\beta$ in $\Homeo(Y)$.
If such a sequence exists, we can construct an asymptotic
homomorphism between the associated crossed product
$C^*$-algebras. In \cite{LM}, however, it was shown that this
simple relation is too weak for Cantor minimal systems. This
happens because there are no consistency in the sequence
$\{\sigma_n\}.$  Moreover, we proved in \cite{M3} that a similar
result holds for minimal dynamical systems on the product of the
Cantor set $X$ and the circle $\T$. To obtain a stronger relation,
one should impose additional conditions on the conjugating maps
$\sigma_n.$  We require that $\sigma_n$ eventually induces the
same map on $K$-theory. Approximate $K$-conjugacy was introduced
in \cite{LM} for Cantor minimal systems, and it was shown
that two minimal systems are approximately $K$-conjugate
if and only if the associated crossed products are isomorphic.
In the present paper, we will show that
two minimal systems on $X\times\T$ associated with
cocycles with values in the rotation group are approximately
$K$-conjugate if and only if
there is an order and unit preserving isomorphism
between the $K$-theory of two associated crossed products
which preserves the images of $K$-theory of $C(X\times\T).$
In fact, in the case that both systems are rigid,
when two systems are approximately $K$-conjugate,
the associated crossed products are isomorphic.
But more is true.
They are $C^*$-strongly approximately flip conjugate
(see Definition \ref{4dcapp} below).
The only difference from the case of Cantor minimal
systems is that we have to control a special projection which does
not come from any clopen subsets of $X$. We call that projection
the generalized Rieffel projection. Difference of two generalized
Rieffel projections will be written by a Bott element associated
with two almost commuting unitaries. A technique developed in
\cite[Section 4]{M3} will be used more carefully in order to fix
the position of that projection in the $K_0$-group.

Non-orientation preserving cases are also discussed.
By using a $\Z_2$-extension,
we can untwist a non-orientation preserving system and
obtain an orientation preserving system.
The associated crossed product $C^*$-algebra turns out
to be isomorphic to the fixed point algebra of
a $\Z_2$-action on the crossed product
associated with the orientation preserving system.
Relation of their $K$-groups is studied.

This paper is organized as follows.
In Section 2, we collect notations and terminologies relevant to
this paper and establish a few elementary facts.
In Section 3, we investigate real rank, stable rank and
comparability of projections of the crossed product $C^*$-algebra.
Section 4 is devoted to the case that cocycles take values
in the rotation group.
In Section 5, when cocycles take values in the rotation group,
we will prove that rigidity implies tracial rank zero.
In Section 6, the generalized Rieffel projection is defined.
In Section 7, we show that isomorphism of $K$-groups implies
approximate $K$-conjugacy
when the systems arise from cocycles
with values in the rotation group.
In Section 8, we deal with non-orientation preserving cases.
In Section 9, we treat examples of various cocycles.
\bigskip

\textbf{Acknowledgement}
The first named author would like to
acknowledge the support from a NSF grant. He would also like to
thank the second author for his effort to make this project
possible.
The second named author was supported by
Grant-in-Aid for Young Scientists (B) of
Japan Society for the Promotion of Science.
He is grateful to Yoshimichi Ueda and Hiroyuki Osaka
for helpful advises.

\section{Preliminaries}

Let $A$ be a unital $C^*$-algebra.
We denote the union of $M_k(A)$ for $k=1,2,\dots$ by $M_\infty(A)$.
The $K_0$-group of $A$ is equipped with the order unit $[1_A]$ and
the positive cone $K_0(A)^+$.
A homomorphism $s:K_0(A)\rightarrow\R$ is called a state
if $s$ carries the order unit to one and the positive cone
to nonnegative real numbers.
We write the set of all states by $S(K_0(A))$ and
call it the state space.
Let $T(A)$ denote the set of all tracial states on $A$.
Endowed with the weak-$*$ topology, $T(A)$ is a compact convex set.
The space of real valued affine continuous functions on $T(A)$ is
written by $\Aff(T(A))$.
We denote by $D$ the natural homomorphism from $K_0(A)$
to $\Aff(T(A))$.
Namely, $D([p])(\tau)$ is equal to $\tau(p)$,
where $p$ is a projection of $M_\infty(A)$ and
$\tau$ is a tracial state on $A$.
We say that the algebra $A$ satisfies Blackadar's second fundamental
comparability question
when the order on projections of $M_\infty(A)$ is determined by traces,
that is, if $p,q\in M_\infty(A)$ are projections and
$\tau(p)<\tau(q)$ for all $\tau\in T(A)$, then $p$ is Murray-von Neumann
equivalent to a subprojection of $q$.

Let $X$ be a compact metrizable space.
Equip $\Homeo(X)$ with the topology of pointwise convergence in norm
on $C(X)$. Thus a sequence $\{\alpha_n\}_{n\in\N}$ in $\Homeo(X)$
converges to $\alpha$, if
\[ \lim_{n\rightarrow\infty}
\sup_{x\in X}\lvert f(\alpha_n^{-1}(x))-f(\alpha^{-1}(x))\rvert=0 \]
for every complex valued continuous function $f\in C(X)$.
This is equivalent to say that
\[ \sup_{x\in X}d(\alpha_n(x),\alpha(x)) \]
tends to zero as $n\rightarrow\infty$, where $d(\cdot,\cdot)$ is a
metric which induces the topology of $X$. When $X$ is the Cantor
set, this is also equivalent to say that, for any clopen subset
$U\subset X$, there exists $N\in\N$ such that
$\alpha_n(U)=\alpha(U)$ for all $n\geq N$.

When $\alpha:X\rightarrow X$ is a homeomorphism
on a compact metrizable space $X$,
we denote the crossed product $C^*$-algebra
arising from the dynamical system $\xa$ by $C^*\xa$.
We use the notation $ufu^*=f\circ\alpha^{-1}$,
where $f$ is a function on $X$ and $u$ is the implementing unitary.
If $X$ is infinite and $\alpha$ is minimal, then $C^*\xa$ is
a simple $C^*$-algebra.
We regard $C(X)$ as a subalgebra of $C^*\xa$.
But when we need to emphasize the embedding,
it will be denoted by $j_{\alpha}:C(X)\rightarrow C^*\xa$.
Let $M_\alpha$ denote the set of $\alpha$-invariant probability
measures on $X$.
If $\alpha$ is free, then there exists a canonical bijection
between $M_\alpha$ and the tracial state space $T(C^*\xa)$.
We may identify these spaces.
Let $\yb$ be another dynamical system.
A continuous surjection $F:X\rightarrow Y$ is called a factor map
if $\beta F=F\alpha$.
It is well-known that $F$ yields an affine continuous surjection
$F_*:M_\alpha\rightarrow M_\beta$ by $F_*(\mu)(E)=\mu(F^{-1}(E))$
for $\mu\in M_\alpha$ and $E\subset Y$.
We also remark that $F$ induces a natural embedding
of $C^*\yb$ into $C^*\xa$.

When $\alpha$ is a minimal homeomorphism on the Cantor set $X$,
we call $\xa$ a Cantor minimal system.
We briefly review results of \cite{GPS}.
The crossed product $C^*$-algebra $C^*\xa$ is
a unital simple AT algebra of real rank zero, and so
it can be classified by its $K$-groups.
The $K_1$-group of $C^*\xa$ is isomorphic to $\Z$,
and $K_0(C^*\xa)$ is unital order isomorphic to
\[ K^0\xa=C(X,\Z)/\{f-f\circ\alpha^{-1}:f\in C(X,\Z)\} \]
equipped with the positive cone
\[ K^0\xa^+=\{[f]:f\in C(X,\Z),f\geq0\} \]
and the order unit $[1_X]$,
where $[f]$ means its equivalence class.
We sometimes write $[f]_\alpha$ to specify $\alpha$.
Throughout this paper, $K_0(C^*\xa)$ will be identified
with $K^0\xa$.
By \cite[Theorem 2.1]{GPS}, $K^0\xa$ is a complete invariant
for strong orbit equivalence of Cantor minimal systems.
The idea of Kakutani-Rohlin partitions will be used repeatedly.
We refer the reader to \cite[Theorem 4.2]{HPS} or
\cite[Section 2]{M3} for details.

We identify the circle with $\T=\R/\Z$ and
write the distance from $t\in\T$ to zero in $\T$
by $\lvert t\rvert$.
Let $R_t$ denote the translation on $\T=\R/\Z$ by $t\in\T$.
Then $\{R_t:t\in\T\}$ forms an abelian subgroup of $\Homeo(\T)$.
We call it the rotation group.
The set of isometric homeomorphisms on $\T$ is written
by $\Isom(\T)$.
Thus,
\[ \Isom(\T)=\{R_t:t\in\T\}\cup\{R_t\lambda:t\in\T\}, \]
where $\lambda\in\Homeo(\T)$ is defined by $\lambda(t)=-t$.
The finite cyclic group of order $m$ is denoted
by $\Z_m\cong\Z/m\Z$ and
may be identified with $\{0,1,\dots,m-1\}$.

Define $o:\Homeo(\T)\rightarrow\Z_2$ by
\[ o(\phi)=\begin{cases}
0 & \phi\text{ is orientation preserving}\\
1 & \phi\text{ is orientation reversing.} \end{cases} \]
Then the map $o(\cdot)$ is a homomorphism.
Let $\Homeo^+(\T)$ denote the set of orientation preserving
homeomorphisms.
Note that $\Homeo^+(\T)\cap\Isom(\T)$ consists of rotations.
\bigskip

In the rest of this section,
we establish notation and some elementary facts
concerning dynamical systems on the product of the Cantor set $X$
and the circle $\T$.

\begin{lem}
Let $\gamma$ be a homeomorphism on $X\times\T$.
Then, there exist $\alpha\in\Homeo(X)$ and a continuous map
$\phi:X\rightarrow\Homeo(\T)$ such that
$\gamma(x,t)=(\alpha(x),\phi_x(t))$
for all $(x,t)\in X\times \T$.
\end{lem}
\begin{proof}
This is obvious because the connected component including $(x,t)$
is $\{x\}\times\T$ and it must be carried to
a connected component by the homeomorphism $\gamma$.
\end{proof}

We denote the homeomorphism of the form in the lemma above
by $\alpha\times\phi$ for short.
When $\phi_x=R_{\xi(x)}$ with a continuous function $\xi:X\rightarrow\T$,
we write $\alpha\times R_\xi$.
Let $F:X\times\T\rightarrow X$ be the projection
onto the first coordinate.
Then we have $F\circ(\alpha\times\phi)=\alpha\circ F$.
Thus, $F$ is a factor map from $(X\times\T,\alpha\times\phi)$ to
$\xa$, and so $F$ induces an affine continuous map
from the set of invariant measures of
$(X\times\T,\alpha\times\phi)$ to that of $\xa$.
Note that if $\alpha\times\phi$ is minimal
then $\alpha$ is also minimal.

We say that $\phi,\psi:X\rightarrow\Homeo(\T)$ are cohomologous,
when there exists a continuous map $\omega:X\rightarrow\Homeo(\T)$
such that $\psi_x\omega_x=\omega_{\alpha(x)}\phi_x$ for all $x\in X$.
If $\phi$ and $\psi$ are cohomologous, it can be easily verified that
$\alpha\times\phi$ and $\alpha\times\psi$ are conjugate.

Let $o(\phi)$ be the composition of $\phi:X\rightarrow\Homeo(\T)$
and $o:\Homeo(\T)\rightarrow\Z_2$, that is, $o(\phi)(x)=o(\phi_x)$.
Under the identification of
\[ C(X,\Z_2)/\{f-f\alpha^{-1}:f\in C(X,\Z_2)\} \]
with $K^0\xa/2K^0\xa$ (see \cite[Lemma 3.5]{M1}),
an element of $K^0\xa/2K^0\xa$ is obtained from $o(\phi)$.
We write it by $[o(\phi)]$ or $[o(\phi)]_\alpha$.

\begin{df}
Let $\xa$ be a Cantor minimal system and
$\phi:X\rightarrow\Homeo(\T)$ be a continuous map.
We say that $\alpha\times\phi$ or $\phi$ is orientation preserving,
if $[o(\phi)]$ is zero in $K^0\xa/2K^0\xa$.
\end{df}

Notice that the concept of `orientation reversing' does not
make sense in this situation and that
$(\alpha\times\phi)^2$ may not be orientation preserving.

\begin{lem}[{\cite[Lemma 4.3]{M3}}]\label{orientation}
Let $\xa$ be a Cantor minimal system.
If $\alpha\times\phi$ is an orientation preserving homeomorphism
on $X\times\T$, then
there exists a continuous map $\psi:X\rightarrow\Homeo^+(\T)$
such that $\phi$ is cohomologous to $\psi$.
\end{lem}
\bigskip

Let $\xa$ be a Cantor minimal system and
let $\phi:X\rightarrow\Homeo(\T)$ be a continuous map.
We would like to compute the $K$-groups
of $C^*(X\times\T,\alpha\times\phi)$.
Let us begin with the orientation preserving case.
When $[o(\phi)]$ is zero in $K^0\xa/2K^0\xa$,
by Lemma \ref{orientation}, we may assume that $o(\phi)(x)$ is zero
for all $x\in X$.
It is evident that $\alpha\times\phi$ induces the action
\[ \alpha^*:f\mapsto f\circ\alpha^{-1} \]
on $K_0(C(X\times\T))\cong C(X,\Z)$,
and that the kernel of $\id-\alpha^*$ is spanned by $[1_X]$
and the cokernel of $\id-\alpha^*$ is $K^0\xa$.
Since $o(\phi)(x)=0$ for all $x\in X$, the induced action
on $K_1(C(X\times\T))\cong C(X,\Z)$ is $\alpha^*$, too.
Consequently we have the following.

\begin{lem}\label{Kofop}
Let $\xa$ be a Cantor minimal system and
let $\alpha\times\phi$ be an orientation preserving homeomorphism
on $X\times\T$.
Then both $K_0(C^*(X\times\T,\alpha\times\phi))$ and
$K_1(C^*(X\times\T,\alpha\times\phi))$ are isomorphic to
$\Z\oplus K^0\xa$.
\end{lem}

Of course the embedding $C^*\xa\subset C^*(X\times\T,\alpha\times\phi)$
induces the embedding of $K^0\xa$ into $K_0(C^*(X\times\T,\alpha\times\phi))$
preserving the unit.
When $\alpha\times\phi$ is minimal, it can be easily verified that
this is really an order embedding.
The order structure of the whole $K_0$-group will be apparent
in later sections.

Next, let us consider the non-orientation preserving case.
Suppose that $[o(\phi)]$ is not zero in $K^0\xa/2K^0\xa$.
Clearly $\alpha\times\phi$ induces the same action on
$K_0(C(X\times\T))\cong C(X,\Z)$ as the orientation preserving case.
But, on the $K_1$-group, the induced action is different.
It is given by
\[ \alpha_\phi^*(f)(x)=
(-1)^{o(\phi)(\alpha^{-1}(x))}f(\alpha^{-1}(x)) \]
for $f\in C(X,\Z)$.
We need to know the kernel and the cokernel of $\id-\alpha_\phi^*$.
Suppose that $f\in C(X,\Z)$ belongs to $\Ker(\id-\alpha_\phi^*)$
and $f\neq0$.
Since $\vert f(x)\rvert=\lvert\alpha_\phi^*(f)(x)\rvert
=\lvert\alpha^*(f)(x)\rvert$ for all $x\in X$,
the minimality of $\alpha$ implies that
$\lvert f(x)\rvert$ is a constant function.
Define $c\in C(X,\Z_2)$ by $f(x)=(-1)^{c(x)}\lvert f(x)\rvert$.
Then $f(x)=(-1)^{o(\phi)(\alpha^{-1}(x))}f(\alpha^{-1}(x))$ yields
$c(x)=o(\phi)(\alpha^{-1}(x))+c(\alpha^{-1}(x))$,
which contradicts $[o(\phi)]\neq0$.
Hence $\id-\alpha_\phi^*$ is injective.
It follows that $K_0(C^*(X\times\T,\alpha\times\phi))\cong
K^0(X,\alpha).$

Let $(X\times\Z_2,\alpha\times o(\phi))$ be the skew product extension
associated with the $\Z_2$-valued cocycle $o(\phi)$.
As $[o(\phi)]\neq0$, $(X\times\Z_2,\alpha\times o(\phi))$ is
a Cantor minimal system (see \cite[Lemma 3.6]{M1}).
By definition, $K^0(X\times\Z_2,\alpha\times o(\phi))$ is isomorphic to
the cokernel of $\id-(\alpha\times o(\phi))^*$ on $C(X\times\Z_2,\Z)$.
Let $\pi$ be the projection from $X\times\Z_2$ onto the first coordinate.
Then $\pi$ is a factor map from $(X\times\Z_2,\alpha\times o(\phi))$
to $\xa$.
It is well-known that $[f]\mapsto[f\circ\pi]$ is
an order embedding between the $K^0$-groups, and so we will regard
$K^0\xa$ as a subgroup of $K^0(X\times\Z_2,\alpha\times o(\phi))$.
It is convenient to introduce a monomorphism $\delta$
from $C(X,\Z)$ to $C(X\times\Z_2,\Z)$
defined by $\delta(f)(x,k)=(-1)^kf(x)$.
Then we can check that
\[ \delta\circ\alpha_\phi^*=(\alpha\times o(\phi))^*\circ\delta. \]
Define $\gamma\in\Homeo(X\times\Z_2)$ by $\gamma(x,k)=(x,k+1)$.
Note that $\gamma$ commutes with $\alpha\times o(\phi)$.
Since $\Ima\delta=\Ima(\id-\gamma^*)$, we have
\begin{align*}
\Coker(\id-\alpha_\phi^*)
&\cong \Ima\delta/\Ima\delta\circ(\id-\alpha_\phi^*) \\
&= \Ima(\id-\gamma^*)/
\Ima((\id-(\alpha\times o(\phi))^*)\circ(\id-\gamma^*)) \\
&\cong C(X\times\Z_2,\Z)/
(\Ker(\id-\gamma^*)+\Ima(\id-(\alpha\times o(\phi))^*)) \\
&\cong K^0(X\times\Z_2,\alpha\times o(\phi))/K^0\xa. \end{align*}
We can summarize the conclusion just obtained as follows.

\begin{lem}\label{Kofnop}
Let $\xa$ be a Cantor minimal system and
let $\alpha\times\phi$ be a non-orientation preserving homeomorphism
on $X\times\T$.
Then we have
\[ K_0(C^*(X\times\T,\alpha\times\phi))\cong K^0\xa \]
and
\[ K_1(C^*(X\times\T,\alpha\times\phi))\cong
\Z\oplus K^0(X\times\Z_2,\alpha\times o(\phi))/K^0\xa, \]
where the $\Z$-summand of the $K_1$-group is
generated by the implementing unitary.
\end{lem}

By the same reason as orientation preserving cases, one sees that
$K_0(C^*(X\times\T,\alpha\times\phi))$ is unital order isomorphic to
$K^0\xa$, when $\alpha\times\phi$ is minimal.

We remark that the torsion subgroup of
\[ K^0(X\times\Z_2,\alpha\times o(\phi))/K^0\xa \]
is isomorphic to $\Z_2$ by \cite[Lemma 4.5]{M3}.
The torsion element is given by
\[ f_0(x,k)=\begin{cases}
1 & o(\phi)(\alpha^{-1}(x))=1\text{ and }k=0 \\
0 & \text{otherwise.} \end{cases} \]

\section{Real rank and stable rank}

Let $\xa$ be a Cantor minimal system and
let $X\ni x\mapsto\phi_x\in\Homeo(\T)$ be a continuous map.
In this section, we would like to compute the real rank and
the stable rank of $C^*(X\times\T,\alpha\times\phi)$.
Let us denote the crossed product $C^*$-algebra by $A$
and the implementing unitary by $u$.
Our crucial tool is a ``large'' subalgebra
of $A$, which is described below.
The proof will be done by applying the argument of \cite{Ph2}
to that subalgebra.

We would like to begin with the definition of the rigidity.

\begin{df}
Let $\alpha\times\phi$ be a homeomorphism on the product of
the Cantor set $X$ and the circle $\T$.
We say that $\alpha\times\phi$ is rigid, if the canonical factor map
from $(X\times\T,\alpha\times\phi)$ to $\xa$ induces
an isomorphism between the sets of invariant probability measures.
\end{df}

\begin{rem}\label{rigidbutnotminimal}
In the definition above, even if $\alpha$ is minimal and
$\alpha\times\phi$ is rigid, $\alpha\times\phi$ may not
be minimal. For example, let $\xa$ be an odometer system and
let $\phi\in\Homeo^+(\T)$ be a Denjoy homeomorphism,
that is, the rotation number of $\phi$ is irrational but
$\phi$ is not conjugate to a rotation.
It is well-known that $\phi$ has a unique invariant
nontrivial closed subset
$Y$, and so $\alpha\times\phi$ has
a nontrivial closed invariant subset $X\times Y$.
Thus $\alpha\times\phi$ is not minimal.
It is also known that $Y$ is the Cantor set and $\phi|Y$ is minimal.
This Cantor minimal system is called a Denjoy system
(see \cite{PSS} for details).
The complement of $Y$ consists of countable disjoint open intervals
and each interval is a wandering set.
Hence, for every $\alpha\times\phi$-invariant probability measure
$\mu$, we have $\mu(X\times Y^c)=0$.
Moreover, as pointed out in \cite[Section 7 (2)]{M1},
the product of $\xa$ and $(Y,\phi|Y)$ is uniquely ergodic.
It follows that $\alpha\times\phi$ is uniquely ergodic.
In particular, it is rigid.
\end{rem}

For $x\in X$,
let $A_x$ be the $C^*$-subalgebra generated by
$C(X\times\T)$ and $uC_0((X\setminus\{x\})\times\T)$.
In \cite[Theorem 3.3]{Pu1}, it was proved that $A_x\cap C^*\xa$ is
a unital AF algebra, where we regard $C^*\xa$ as
a $C^*$-subalgebra of $A$.
This AF subalgebra played a crucial role in the subsequent papers
\cite{HPS} and \cite{GPS}.
In our situation, the $C^*$-subalgebra $A_x$ is not an AF algebra
but an AT algebra,
and it helps us to show real rank zero.

\begin{prop}\label{Putcutdown}
In the setting above, we have the following.
\begin{enumerate}
\item $A_x$ is a unital AT algebra.
\item When $\alpha\times\phi$ is orientation preserving,
$K_0(A_x)$ is unital order isomorphic to $K^0\xa$ and
$K_1(A_x)$ is isomorphic to $K^0\xa$.
\item When $\alpha\times\phi$ is not orientation preserving,
$K_0(A_x)$ is unital order isomorphic to $K^0\xa$ and
$K_1(A_x)$ is isomorphic to an extension of
$\Coker(\id-\alpha_\phi^*)$ by $\Z$.
\item There exists a canonical bijection between
the tracial state space $T(A_x)$ and
the set of $\alpha\times\phi$-invariant probability measures.
\item $A_x$ is simple if and only if $\alpha\times\phi$ is minimal.
\item If $\alpha\times\phi$ is minimal, then
$A_x$ is real rank zero if and only if
$\alpha\times\phi$ is rigid.
\end{enumerate}
\end{prop}
\begin{proof}
(1) Let
\[ {\cal P}_n=\{X(n,v,k):v\in V_n,k=1,2,\dots,h(v)\} \]
be a sequence of Kakutani-Rohlin partitions which gives
a Bratteli-Vershik model for $\xa$ (see \cite[Theorem 4.2]{HPS} or
\cite[Section 2]{M3} for Kakutani-Rohlin partitions).
We assume that the sequence of the roof sets
\[ R({\cal P}_n)=\bigcup_{v\in V}X(n,v,h(v)) \]
shrinks to a singleton $\{x\}$.
Let $A_n$ be the $C^*$-subalgebra generated by
$C(X\times\T)$ and $uC(R({\cal P}_n)^c\times\T)$.
It is easy to see that
$A_x$ is the norm closure of the union of all $A_n$'s.
By using a similar argument to \cite[Lemma 3.1]{Pu1},
it can be shown that $A_n$ is isomorphic to
\[ \bigoplus_{v\in V_n}
M_{h(v)}\otimes C(X(n,v,h(v)))\otimes C(\T), \]
which is an AT algebra.
To verify this, define a projection $p_v$ by
\[ p_v=\sum_{k=1}^{h(v)}1_{X(n,v,k)} \]
for each $v\in V$.
Since $p_vu(1-1_{R({\cal P}_n)})=u(1-1_{R({\cal P}_n)})p_v$,
the projection $p_v$ is central in $A_n$.
Clearly $u^{i-j}1_{X(n,v,j)}$ for $i,j=1,2,\dots,h(v)$ form matrix units
of $p_vA_np_v$.
By
\[ 1_{X(n,v,h(v))}A_n1_{X(n,v,h(v))}=C(X(n,v,h(v))\times\T), \]
we obtain the description above.
Hence $A_x$ is also an AT algebra.

(2) There is a natural homomorphism from $K_i(C(X\times\T))\cong C(X,\Z)$
to $K_i(A_n)$ for $i=1,2$.
By (1), the kernel of this map is
\[ \{f-f\circ\alpha^{-1}:f\in C(X,\Z),
f(y)=0\text{ for all }y\in R({\cal P}_n)\}. \]
Therefore $K_i(A_x)$ is isomorphic to
\[ C(X,\Z)/\{f-f\circ\alpha^{-1}:f\in C(X,\Z),f(x)=0\}. \]
It follows from $1_X\circ\alpha^{-1}=1_X$ that
\[ \{f-f\circ\alpha^{-1}:f\in C(X,\Z),f(x)=0\}=
\{f-f\circ\alpha^{-1}:f\in C(X,\Z)\}, \]
which implies $K_i(A_x)\cong K^0\xa$.

(3) The computation of $K_0$-group is the same as
the orientation preserving case.
Let us consider $K_1(A_x)$.
It is not hard to see that $K_1(A_x)$ is isomorphic to
\[ C(X,\Z)/
\{f-\alpha_\phi^*(f):f\in C(X,\Z)\text{ and }f(x)=0\}. \]
We follow the notation used in the discussion
before Lemma \ref{Kofnop}.
The image of
\[ \{f-\alpha_\phi^*(f):f\in C(X,\Z)\text{ and }f(x)=0\} \]
by $\delta$ is equal to the image of
\[ \{f-f\circ(\alpha\times o(\phi))^{-1}:
f\in C(X\times\Z_2,\Z),f(x,0)=(x,1)=0\} \]
by $\id-\gamma^*$. Hence, in the same way as Lemma \ref{Kofnop},
we have
\[ K_1(A_x)\cong
K^0(X\times\Z_2,\alpha\times o(\phi);\Z_2)/K^0\xa. \]
See \cite{M2} for the definition of
$K^0(X\times\Z_2,\alpha\times o(\phi);\Z_2)$.
By \cite[Theorem 4.1]{Pu1},
\[ 0\rightarrow\Z\rightarrow
K^0(X\times\Z_2,\alpha\times o(\phi);\Z_2)\rightarrow
K^0(X\times\Z_2,\alpha\times o(\phi))\rightarrow 0 \]
is exact, which implies that $K_1(A_x)$ is an extension of
\[ \Coker(\id-\alpha_\phi^*)\cong
K^0(X\times\Z_2,\alpha\times o(\phi))/K^0\xa \]
by the integers $\Z$.

(4) Since there exists a canonical bijection between $T(A)$ and
the set of $\alpha\times\phi$-invariant probability measures,
it suffices to check that every $\tau\in T(A_x)$ extends to
a tracial state on $A$.
Let $f\in C(X\times\T)$ be a function satisfying $0\leq f\leq1$.
Take a natural number $n\in\N$ arbitrarily.
We can find a clopen neighborhood $U$ of $x$ such that
$U,\alpha^{-1}(U),\dots,\alpha^{-n}(U)$ are mutually disjoint.
Let $p=1_{U\times\T}\in C(X\times\T)$.
Then we have $\tau(p)<n^{-1}$,
because $p,u^*pu,\dots,u^{n*}pu^n$ are mutually equivalent in $A_x$
and mutually disjoint.
We also notice that $u(1-p)f$ belongs to $A_x$.
It follows that
\[ \tau(f)\leq\tau(p)+\tau((1-p)f)=\tau(p)+\tau(u(1-p)fu^*)
<\frac{1}{n}+\tau(ufu^*). \]
Similarly we can see $\tau(ufu^*)<n^{-1}+\tau(f)$.
Since $n$ is arbitrary, $\tau(f)$ equals $\tau(ufu^*)$,
which means that $\tau$ extends to a trace on $A$.

(5) Note that the $C^*$-algebra $A$ can be regarded as
a groupoid $C^*$-algebra associated with the equivalence relation
\[ {\cal R}=\{(z,(\alpha\times\phi)^k(z)):
z\in X\times\T,k\in\Z\}. \]
Then the $C^*$-subalgebra $A_x$ corresponds to the subequivalence
relation
\[ {\cal R}_x={\cal R}\setminus
\{((\alpha\times\phi)^k(x,t),(\alpha\times\phi)^l(x,t)):
t\in\T, (1-k,l)\in\N^2\text{ or }(k,1-l)\in\N^2\}. \]
It is well-known that a homeomorphism on a compact space is minimal
if and only if every positive orbit is dense.
Therefore $\alpha\times\phi$ is minimal if and only if
each equivalence class of ${\cal R}_x$ is dense in $X\times\T$.
It follows from \cite[Proposition 4.6]{R} that this is equivalent to
$A_x$ being simple.

(6) By (1) and (5), $A_x$ is a unital simple AT algebra.
By (2) and (3), in the $K_0$-group,
every projection of $A_x$ is equivalent to some $[f]\in K^0\xa$.
Hence, projections in $A_x$ separate traces on $A_x$
if and only if $\alpha\times\phi$ is rigid.
Then the conclusion follows from \cite[Theorem 1.3]{BBEK}.
\end{proof}

\begin{rem}
The six-term exact sequence of \cite[Theorem 2.4]{Pu2}
applies to this situation.
See \cite[Example 2.6]{Pu2}.
The reduced groupoid $C^*$-algebra $C^*_r(H)$ appearing there
is isomorphic to $C(\T)\otimes{\cal K}$.
\end{rem}
\bigskip

We would like to consider the real rank of $A$.
In \cite{Ph2}, it was shown that if $G$ is an almost AF Cantor
groupoid and $C_r^*(G)$ is simple,
then $C_r^*(G)$ has real rank zero.
The key of its proof was the presence of a ``large'' AF subalgebra.
We will show that a similar argument is possible
when the AF subalgebra is replaced
by a subalgebra with tracial rank zero.

\begin{df}[{\cite[Theorem 6.13]{L3}}]
We say that a unital simple $C^*$-algebra $A$ has
tracial (topological) rank zero,
if for any finite subset ${\cal F}\subset A$, any $\varepsilon>0$
and any nonzero positive element $c\in A$,
there exists a projection $e\in A$ and a finite dimensional
unital subalgebra $E\subset eAe$
(that is, $e$ is the identity of $E$) such that:
\begin{enumerate}
\item $\lVert ae-ea\rVert<\varepsilon$ for all $a\in{\cal F}$.
\item For every $a\in{\cal F}$, there is $b\in E$ such that
$\lVert pap-b\rVert<\varepsilon$.
\item $1-e$ is Murray-von Neumann equivalent to a projection
in $\overline{cAc}$.
\end{enumerate}
\end{df}

\begin{df}\label{dfofalmostTAF}
A unital simple $C^*$-algebra $A$ is called
an almost tracially AF algebra,
if there exists a unital simple subalgebra $B\subset A$
with tracial rank zero
such that the following holds:
for any finite subset ${\cal F}\subset A$ and any $\varepsilon>0$,
there exists a projection $p\in B$ such that:
\begin{enumerate}
\item For every $a\in {\cal F}$, there is $b\in B$ such that
$\lVert ap-b\rVert<\varepsilon$.
\item $\tau(1-p)<\varepsilon$ for every tracial state $\tau\in T(B)$.
\end{enumerate}
\end{df}

\begin{lem}\label{almostTAF}
Suppose that $\alpha\times\phi$ is a minimal homeomorphism
on $X\times\T$.
If $\alpha\times\phi$ is rigid, then
the crossed product $C^*$-algebra
$A=C^*(X\times\T,\alpha\times\phi)$ is an almost tracially AF algebra.
\end{lem}
\begin{proof}
Take $x\in X$. By Proposition \ref{Putcutdown},
$A_x$ is a unital simple AT algebra of real rank zero.
It follows from \cite[Proposition 2.6]{L1} that $A_x$ is
a unital simple algebra with tracial rank zero.
Suppose that a finite subset ${\cal F}\subset A$ and
$\varepsilon>0$ are given.
Since
\[ \left\{\sum_{k=-N}^Nf_ku^k:N\in\N, f_k\in C(X\times\T)\right\} \]
is a dense subalgebra of $A$,
we may assume that there is $N\in\N$ such that
${\cal F}$ is contained in
\[ E_N=\left\{\sum_{k=-N}^Nf_ku^k:f_k\in C(X\times\T)\right\}. \]
We can find a clopen neighborhood $U$ of $x$ so that
$\alpha^{1-N}(U),\alpha^{2-N}(U),\dots,U,\alpha(U),\dots,\alpha^N(U)$
are mutually disjoint and
$\mu(U)<\varepsilon/2N$ for all $\mu\in M_\alpha$.
Put $V=\bigcup_{k=1-N}^N\alpha^k(U)$ and
$p=1_{V^c\times\T}\in C(X\times\T)$.
It is easy to see that $u^kp$ and $u^{1-k}p$ belong to $A_x$
for $k=0,1,\dots,N$,
and so $ap\in A_x$ for every $a\in E_N$.
Moreover, $\tau(1-p)$ is less than $\varepsilon$
for every $\tau\in T(A_x)$.
This finishes the proof of the lemma.
\end{proof}

The following lemma is a generalization of \cite[Lemma 4.3]{Ph2}.

\begin{lem}\label{keyofrr0}
Let $A$ be a unital simple algebra with tracial rank zero.
Let $p\in A$ be a projection and
let $a\in A$ be a nonzero self-adjoint element.
For any $\varepsilon>0$ and $n\in\N$ satisfying $n\varepsilon>1$,
there exists a projection $q\in A$ such that
$\lVert qa-aq\rVert<\varepsilon\lVert a\rVert$, $p\leq q$ and
\[ \tau(q)<(2n+1)\tau(p) \]
for all $\tau\in T(A)$.
\end{lem}
\begin{proof}
Without loss of generality, we may assume $\lVert a\rVert=1$.
Put $\varepsilon_0=8^{-1}(\varepsilon-n^{-1})$.
Choose $\delta_0>0$ so that whenever $p,q\in A$ are projections
satisfying $\lVert pq-p\rVert<\delta_0$,
then there exists a unitary $u\in A$ such that
\[ \lVert u-1\rVert<\varepsilon_0 \ \text{ and } \ p\leq uqu^*. \]
Let
\[ \delta=\min\{\varepsilon_0,4^{-1}\delta_0,
(2n+1)^{-1}\tau(p):\tau\in T(A)\}. \]
By \cite[Proposition 2.4]{L1} together with
results of \cite[Section 3]{L1},
there exist a projection $e\in A$ and
a finite dimensional unital subalgebra $E\subset eAe$ such that:
\begin{itemize}
\item $\lVert ae-ea\rVert<\delta$ and $\lVert pe-ep\rVert<\delta$.
\item There exist $b,c\in E$ such that $\lVert eae-b\rVert<\delta$
and $\lVert epe-c\rVert<\delta$.
\item $\tau(1-e)<\delta$ for every tracial state $\tau\in T(A)$.
\end{itemize}
We may assume that $b=b^*$, $\lVert b\rVert=1$ and $c$ is a projection.
Thanks to \cite[Lemma 4.2]{Ph2}, there exists a projection $q_0\in E$
such that
\[ c\leq q_0, \ [q_0]\leq 2n[c]\in K_0(E), \
\text{ and } \lVert q_0b-bq_0\rVert<\frac{1}{n}. \]
Put $q=1-e+q_0$.
Since
\begin{align*}
& \lVert pq-p\rVert \\
&\leq \lVert(p(1-e)+peq_0)-(p(1-e)+pcq_0)\rVert+
\lVert(p(1-e)+pcq_0)-(p(1-e)+pe)\rVert \\
&< 2\delta+2\delta\leq\delta_0,
\end{align*}
there is a unitary $u\in A$ such that $\lVert u-1\rVert<\varepsilon_0$
and $p\leq uqu^*$.
It is not hard to see that
\begin{align*}
& \tau(uqu^*)=\tau(q)=\tau(1-e)+\tau(q_0) \\
&< \delta+2n\tau(c)<(2n+1)\delta+2n\tau(pe)<(2n+1)\tau(p)
\end{align*}
for all $\tau\in T(A)$ and that
\begin{align*}
\lVert [uqu^*,a]\rVert&< \lVert[q,a]\rVert+4\varepsilon_0 \\
&< \lVert[q,b]\rVert+4\delta+4\varepsilon_0 \\
&< \frac{1}{n}+8\varepsilon_0=\varepsilon,
\end{align*}
thereby completing the proof.
\end{proof}

The following is a well-known matrix trick.
We omit the proof.

\begin{lem}\label{trick}
When $a$ is a self-adjoint element of a unital $C^*$-algebra $A$,
\[ \begin{bmatrix}a&0\\0&0\end{bmatrix} \]
is approximated by an invertible self-adjoint element
in $A\otimes M_2$.
\end{lem}

Although the proof of the following theorem is almost
the same as that of \cite[Theorem 4.7]{Ph2},
we would like to state it for the reader's convenience.

\begin{thm}
If a unital simple $C^*$-algebra $A$ is an almost tracially AF algebra,
then $A$ has real rank zero.
\end{thm}
\begin{proof}
Let $B\subset A$ be a unital simple subalgebra with tracial rank zero
as in Definition \ref{dfofalmostTAF}.

Let $a\in A$ be self-adjoint and non-invertible.
It suffices to show that $a$ is approximated
by a self-adjoint invertible element of $A$.
Without loss of generality, we may assume $\lVert a\rVert\leq1$.
Take $\varepsilon>0$ arbitrarily.
Define a continuous function $g$ on $[-1,1]$ by
\[ g(t)=\begin{cases}
1-\varepsilon^{-1}\lvert t\rvert & \lvert t\rvert\leq\varepsilon \\
0 & \text{otherwise}. \end{cases} \]
Put
\[ \varepsilon_0=\min\{\tau(g(a)):\tau\in T(A)\}. \]
Since $A$ is simple, $\varepsilon_0$ is positive.
Applying \cite[Lemma 4.4]{Ph2} to $g:[-1,1]\rightarrow[0,1]$ and
$4^{-1}\varepsilon_0>0$, we obtain $\delta>0$.
We may assume that $\delta$ is less than $\varepsilon$.
Choose a natural number $n\in\N$ so that
\[ \frac{1}{n}<\min\left\{\frac{\varepsilon_0}{12},
\frac{\delta}{2}\right\}. \]
By definition, there is a projection $p\in B$ such that
$a(1-p)$ is close to $B$ and $\tau(p)$ is less than $n^{-2}$
for all $\tau\in T(B)$.
By perturbing $a$, we may assume that $a(1-p)$ belongs to $B$.
Then, we can apply Lemma \ref{keyofrr0} to $a-pap\in B$ and
$p\in B$, and get a projection $q\in B$ such that $p\leq q$,
\[ \lVert[q,a-pap]\rVert<
\frac{\delta}{2}\lVert a-pap\rVert\leq\delta \]
and
\[ \tau(q)\leq(2n+1)\tau(p)<\frac{3}{n}
<\frac{\varepsilon_0}{4} \]
for every tracial state $\tau\in T(B)$.
It follows that
\[ \tau(g(a))-\tau(q)>\varepsilon_0-\frac{\varepsilon_0}{4}
>\frac{\varepsilon_0}{4}. \]
We also notice that $\lVert[q,a]\rVert<\delta$.
Put $a_0=(1-q)a(1-q)\in B$.
By the choice of $\delta$, we have
\[ \tau(g(a_0))>\tau(g(a))-\tau(q)-\frac{\varepsilon_0}{4}
>\frac{\varepsilon_0}{2}. \]
A unital simple algebra with tracial rank zero has real rank zero
by \cite[Theorem 3.4]{L1}, and so
there exists a projection $r$ in the hereditary subalgebra
of $B$ generated by $g(a_0)$ such that
\[ \lVert rg(a_0)-g(a_0)\rVert<\frac{\varepsilon_0}{4}. \]
Then
\[ \tau(r)\geq\tau(rg(a_0)r)
>\tau(g(a_0))-\frac{\varepsilon_0}{4}
>\frac{\varepsilon_0}{4} \]
for all $\tau\in T(B)$.
Since the order on projections of $B$ is determined by traces
(see \cite[Theorem 6.8, 6.13]{L2}),
there is a projection $r_0\in B$ such that $r_0\leq r$ and
$r_0\sim q$.
Moreover, by means of \cite[Lemma 4.6]{Ph2}, we have
\[ \lVert r_0a_0-a_0r_0\rVert<2\varepsilon \ \text{ and } \
\lVert r_0a_0r_0\rVert<\varepsilon. \]
As a result,
\begin{align*}
& a\stackrel{2\varepsilon}{\approx}qaq+a_0 \\
&\stackrel{4\varepsilon}{\approx}qaq+r_0a_0r_0+
(1-q-r_0)a_0(1-q-r_0) \\
&\stackrel{\varepsilon}{\approx}qaq+(1-q-r_0)a_0(1-q-r_0)
\end{align*}
is obtained.
The element $(1-q-r_0)a_0(1-q-r_0)$ belongs to $B$ and
$B$ has real rank zero.
By applying Lemma \ref{trick} to $qaq$ and $r_0\sim q$,
we can get the conclusion.
\end{proof}

\begin{cor}\label{rr0}
For a minimal homeomorphism $\alpha\times\phi$
on $X\times\T$, the following are equivalent.
\begin{enumerate}
\item $\alpha\times\phi$ is rigid.
\item The crossed product $C^*$-algebra
$A=C^*(X\times\T,\alpha\times\phi)$ has
real rank zero.
\item $D(K_0(A))$ is uniformly dense in $\Aff(T(A))$.
\end{enumerate}
\end{cor}
\begin{proof}
(1)$\Rightarrow$(2).
This is immediate from Lemma \ref{almostTAF} and
the theorem above.

(2)$\Rightarrow$(3).
Since $\alpha\times\phi$ is minimal, $A$ is simple.
By Theorem \ref{sr1}, $A$ is stably finite and
the projections in $A\otimes{\cal K}$ satisfy
cancellation.
Furthermore, $K_0(A)$ is weakly unperforated
by \cite[Theorem 4.5]{Ph1}.
It follows from \cite[Theorem 6.9.3]{B} that
the image of $K_0(A)$ is uniformly dense in
real valued affine continuous functions on $QT(A)$,
the set of quasitraces on $A$.
Because every element of $\Aff(T(A))$ comes from
a self-adjoint element of $A$
(see \cite[Proposition 3.12]{BKR} for instance),
it extends to a real valued affine continuous function on $QT(A)$.
Hence $D(K_0(A))$ is uniformly dense in $\Aff(T(A))$.
Note that if one uses the deep result obtained by Haagerup
in \cite{H}, the latter half of the proof is superfluous.

(3)$\Rightarrow$(1).
Suppose that $\alpha\times\phi$ is not rigid.
Thus, there are $\nu_1,\nu_2\in M_{\alpha\times\phi}$ such that
$F_*(\nu_1)=F_*(\nu_2)=\nu\in M_\alpha$,
where $F$ is the canonical factor map onto $\xa$.
Let $\tau_1$ and $\tau_2$ be the tracial states on $A$
arising from $\nu_1$ and $\nu_2$.
Since
\[ \tau_{1*}([f])=\nu_1(f\circ F)=\nu(f)
=\nu_2(f\circ F)=\tau_{2*}([f]) \]
for all $[f]\in K^0\xa$, projections in $C^*\xa$ cannot separate
traces on $A$.
Therefore, we can finish the proof here
when $\alpha\times\phi$ is not orientation preserving.

Assume that $\alpha\times\phi$ is orientation preserving.
By Lemma \ref{Kofop}, $K_0(A)$ is isomorphic to $\Z\oplus K^0\xa$.
Suppose that
there are projections $e_1,e_2\in M_k(A)$ for some integer $k$ such that
$[e_1]-[e_2]=(1,0)\in\Z\oplus K^0\xa.$
Then, for any $x=(n,[f])\in K_0(A)\cong\Z\oplus K^0\xa$, we have
\[ \tau_{1*}(x)-\tau_{2*}(x)=n(\tau_1(e_1-e_2)-\tau_2(e_1-e_2)). \]
Hence
\[ \{\tau_{1*}(x)-\tau_{2*}(x):x\in K_0(A)\} \]
is discrete in $\R$, which is a contradiction.
\end{proof}
\bigskip

We now turn to a consideration of stable rank of $A$.
Suppose that $\alpha\times\phi$ is minimal but may not be rigid.
The $C^*$-subalgebra $A_x$ may not have real rank zero.
But, it is still a unital simple AT algebra
by Proposition \ref{Putcutdown}.
Moreover, Lemma \ref{almostTAF} is also valid
when one replaces `almost tracially AF' by `almost AT'.
A unital simple AT algebra is known to have property (SP),
that is, every nonzero hereditary subalgebra contains
a nonzero projection.
Hence we see that $A$ also has property (SP)
by virtue of Lemma \ref{almostTAF}.
We also remark that $A_x$ has stable rank one and
the order on projections of $A_x$ is determined by traces,
because $A_x$ is a unital simple AT algebra.
Then, by reading the proof of \cite[Theorem 5.2]{Ph2}
carefully, it turns out that
$A$ and the ``large'' subalgebra $A_x$
do not need to have real rank zero and
that they only need to have property (SP)
so that the proof works.
As a consequence, we have the following.

\begin{thm}\label{sr1}
When $\alpha\times\phi$ is a minimal homeomorphism
on $X\times\T$, the crossed product $C^*$-algebra
$A=C^*(X\times\T,\alpha\times\phi)$ has stable rank one.
In particular, the projections in $A\otimes{\cal K}$ satisfy
cancellation.
\end{thm}

In \cite[Theorem 4.5]{Ph1}, it was proved that
$A$ satisfies the $K$-theoretic version of Blackadar's second
fundamental comparability question, that is,
if $x\in K_0(A)$ satisfies $\tau_*(x)>0$
for all $\tau\in T(A)$, then $x\in K_0(A)^+$.
In particular, $K_0(A)$ is weakly unperforated.
Combining this with the theorem above,
we can deduce the following.

\begin{thm}\label{BSFC}
When $\alpha\times\phi$ is a minimal homeomorphism
on $X\times\T$, the order on projections of
$M_\infty(C^*(X\times\T,\alpha\times\phi))$ is determined by traces.
In other words, $C^*(X\times\T,\alpha\times\phi)$ satisfies
Blackadar's second fundamental comparability question.
\end{thm}

\section{Cocycles with values in the rotation group}

Let $\xa$ be a Cantor minimal system and
let $\xi:X\rightarrow\T$ be a continuous map.
In this section, we would like to investigate a homeomorphism
$\alpha\times R_\xi$ on $X\times\T$ and
its related crossed product $C^*$-algebra
$A=C^*(X\times\T,\alpha\times R_\xi)$.
Of course, $\alpha\times R_\xi$ is orientation preserving.

\begin{df}\label{KT}
Let $\alpha$ be a minimal homeomorphism on $X$.
Define
\[ K^0_\T\xa \ = \ C(X,\T)/
\{\eta-\eta\alpha^{-1}:\eta\in C(X,\T)\}. \]
The equivalence class of $\xi\in C(X,\T)$ in
$K^0_\T\xa$ is denoted by $[\xi]_\alpha$ or $[\xi]$.
\end{df}

Let $\theta\in\T$ and put $\xi(x)=\theta$ for all $x\in X$.
Thus, $\xi$ is a constant function.
It is easy to see that $[\xi]$ is zero in $K^0_\T\xa$
if and only if $\theta$ is a topological eigenvalue
of $\xa$.
The reader may refer to \cite[Theorem 5.17]{W}
for topological eigenvalues.
Since $X$ is compact, the set of topological eigenvalues
is at most countable.
It follows that $K^0_\T\xa$ is uncountable.

At first, we describe when $\alpha\times R_\xi$ is minimal
in terms of $K^0_\T\xa$.
Note that more general results were obtained in \cite{Pa}.

\begin{lem}\label{opminimal}
Let $\xa$ be a Cantor minimal system and
$\xi:X\rightarrow\T$ be a continuous map.
Then, $\alpha\times R_\xi$ is minimal if and only if
$n[\xi]\neq 0$ in $K^0_\T\xa$ for all $n\in\N$.
\end{lem}
\begin{proof}
Suppose that there exist $n\in\N$ and $\eta$ such that
$n\xi=\eta-\eta\alpha^{-1}$.
Then
\[ \{(x,t)\in X\times\T:nt=\eta(\alpha^{-1}(x))\} \]
is closed and invariant under $\alpha\times R_\xi$,
and so $\alpha\times R_\xi$ is not minimal.

Let us prove the other implication.
Assume that $\alpha\times R_\xi$ is not minimal.
Let $E$ be a minimal subset of $\alpha\times R_\xi$.
Note that $\id\times R_t$ commutes with $\alpha\times R_\xi$.
Since $E$ is not the whole of $X\times\T$,
\[ G=\{t\in\T:(\id\times R_t)(E)=E\} \]
is a closed proper subgroup of $\T$.
It follows that there exists $n\in\N$ such that $G=\{t\in\T:nt=0\}$.
Moreover, on account of the minimality of $E$,
we can deduce that there exists $\eta:X\rightarrow\T$ such that
\[ E=\{(x,t)\in X\times\T:nt=\eta(x)\}. \]
The map $\eta$ is continuous, because $E$ is closed.
Hence we have $n\xi=\eta\alpha-\eta$.
\end{proof}

If $\alpha\times R_\xi$ is not minimal,
then there exist uncountably many minimal closed subsets.
In particular, it is not rigid.
Compare this with Remark \ref{rigidbutnotminimal}.

\begin{lem}\label{conju}
Let $\xa$ and $\yb$ be Cantor minimal systems.
Let $\xi$ and $\zeta$ be continuous maps from $X$ to
$\T$. Suppose that $\alpha\times R_\xi$ is minimal.
Then, $\alpha\times\xi$ and $\beta\times\zeta$ is conjugate
if and only if there exists a homeomorphism $F:X\rightarrow Y$
such that $F\alpha=\beta F$, and $[\xi]_\alpha=[\zeta F]_\alpha$
or $[\xi]_\alpha=-[\zeta F]_\alpha$ in $K^0_\T\xa$.
\end{lem}
\begin{proof}
The `if' part is clear. We consider the `only if' part.
Let $F\times\phi:X\times\T\rightarrow Y\times\T$ be
a conjugating map, that is,
\[ F\alpha=\beta F \ \text{ and } \
\phi_{\alpha(x)}(s+\xi(x))=\phi_x(s)+\zeta(F(x)) \]
for all $(x,s)\in X\times\T$.
For every $t\in\T$, $\id\times R_t$ commutes with $\beta\times R_\zeta$,
and so $(F\times\phi)^{-1}(\id\times R_t)(F\times\phi)$ commutes with
$\alpha\times R_\xi$.
Let $x\in X$ and put $s=\phi_x^{-1}(\phi_x(0)+t)$.
Then we have
\[ (F\times\phi)^{-1}(\id\times R_t)(F\times\phi)(x,0)=
(x,s)=(\id\times R_s)(x,0). \]
By the minimality of $\alpha\times\xi$, we can conclude that
\[ (F\times\phi)^{-1}(\id\times R_t)(F\times\phi)=\id\times R_s. \]
It follows that the mapping $t\mapsto s$ is
a continuous injective homomorphism from $\T$ to $\T$.
Thus,
\[ (F\times\phi)^{-1}(\id\times R_t)(F\times\phi)=\id\times R_t \]
for all $t\in\T$, or
\[ (F\times\phi)^{-1}(\id\times R_t)(F\times\phi)=\id\times R_{-t} \]
for all $t\in\T$.
Without loss of generality we may assume the first,
which yields
\[ \phi_x(s+t)=\phi_x(s)+t \]
for all $(x,s)\in X\times\T$ and $t\in\T$.
It follows that $\phi_x$ equals $R_{\phi_x(0)}$ and
\[ \xi(x)+\phi_{\alpha(x)}(0)=\phi_x(0)+\zeta(F(x)). \]
Thereby the assertion follows.
\end{proof}

Next, we would like to consider when $\alpha\times R_\xi$ is rigid.
Although the following is a special case of \cite[Theorem 3]{Pa}
or \cite[Theorem 3.5]{Z},
we include the proof for the reader's convenience.

\begin{lem}\label{whenrigid}
Let $\xa$ be a Cantor minimal system.
For a continuous function $\xi:X\rightarrow\T$,
the following are equivalent.
\begin{enumerate}
\item $\alpha\times R_\xi$ is rigid.
\item Every $\alpha\times R_\xi$-invariant measure $\nu$
is also $\id\times R_t$-invariant for all $t\in\T$, that is,
$\nu$ is a product measure of the Haar measure on $\T$
and an invariant measure for $\xa$.
\item For every $\alpha$-invariant measure $\mu$
on $X$ and $n\in\N$,
there does not exist a Borel function
$\eta:X\rightarrow\T$ such that
$n\xi(x)=\eta(x)-\eta\alpha^{-1}(x)$ for $\mu$-almost every $x\in X$.
\end{enumerate}
\end{lem}
\begin{proof}
We denote the canonical factor map from $(X\times\T,\alpha\times R_\xi)$
to $\xa$ by $\pi$.

(1)$\Rightarrow$(2).
This is immediate from
$\pi_*\circ(\id\times R_t)_*=(\pi\circ(\id\times R_t))_*=\pi_*$.

(2)$\Rightarrow$(1).
Define a continuous map $\Phi:C(X\times\T)\rightarrow C(X)$ by
\[ \Phi(f)(x)=\int_\T f(x,t)\,dt. \]
If $\nu\in M_{\alpha\times R_\xi}$ is $\id\times R_t$-invariant,
then $\nu(f)=\pi_*(\nu)(\Phi(f))$ for all $f\in C(X\times\T)$.
Hence $\pi_*^{-1}(\pi_*(\nu))=\{\nu\}$.

(2)$\Rightarrow$(3).
Suppose that there exist an $\alpha$-invariant measure
$\mu\in M_\alpha$, $n\in\N$ and a Borel function $\eta:X\rightarrow\T$
such that
\[ n\xi(x)=\eta(x)-\eta\alpha^{-1}(x) \]
for $\mu$-almost every $x\in X$. Then
\[ C(X\times\T)\ni f\mapsto
\frac{1}{n}\int_X\sum_{nt=\eta\alpha^{-1}(x)}f(x,t)\,d\mu(x)\in\C \]
yields a probability measure on $X\times\T$.
Note that the summation runs over $n$ distinct $t$'s
which satisfy $nt=\eta\alpha^{-1}(x)$.
This measure is $\alpha\times R_\xi$-invariant,
because we have
\[ \sum_{nt=\eta\alpha^{-1}(x)}f(\alpha(x),t+\xi(x))
=\sum_{nt=\eta\alpha^{-1}(\alpha(x))}f(\alpha(x),t) \]
for $\mu$-almost every $x\in X$.
But it is not the product of the Haar measure and $\mu$.

(3)$\Rightarrow$(2).
Suppose that $\nu\in M_{\alpha\times R_\xi}$ is not invariant
under the rotation $\id\times R_t$.
We may assume that $\nu$ is an ergodic measure.
It follows from \cite[Corollary 4.1.9]{KH} or \cite[Lemma 6.13]{W}
that there exists
an $F_\sigma$ subset $E\subset X\times\T$ such that
$\nu(E)=1$ and
\[ \lim_{n\rightarrow\infty}\frac{1}{n}\sum_{k=0}^{n-1}
f((\alpha\times R_\xi)^k(x,s))=\nu(f) \]
for all $f\in C(X\times\T)$ and $(x,s)\in E$.
We may assume that $E$ is $\alpha\times R_\xi$-invariant.
Put
\[ G=\{t\in\T:(\id\times R_t)_*(\nu)=\nu\}. \]
By assumption, $G$ is a closed proper subgroup of $\T$.
Thus, there is $n\in\N$ such that $G=\{t\in\T:nt=0\}$.
If $(x,s)$ belongs to $E$, then we have
\[ \lim_{n\rightarrow\infty}\frac{1}{n}\sum_{k=0}^{n-1}
f((\alpha\times R_\xi)^k(x,s+t))=\nu(f\circ(\id\times R_t))
=(\id\times R_t)_*(\nu)(f) \]
for all $f\in C(X\times\T)$.
Therefore $(\id\times R_t)(E)\cap E$ is empty if $t\notin G$.
Furthermore, by replacing $E$ by
\[ \bigcup_{t\in G}(\id\times R_t)(E), \]
we may assume that $(\id\times R_t)(E)=E$ for all $t\in G$.
On the $F_\sigma$ subset $\pi(E)\subset X$,
we define a $\T$-valued function $\eta$ by
$\eta(x)=nt$ for $(x,t)\in E$.
If $E_0\subset E$ is closed, then
$\eta$ is evidently continuous on $\pi(E_0)$.
It follows that $\eta$ is a well-defined Borel function.
For $(x,s)\in E$,
\[ (\alpha\times R_\xi)(x,s)=(\alpha(x),s+\xi(x)) \]
belongs to $E$, and so
\[ \eta(\alpha(x))=ns+n\xi(x)=\alpha(x)+n\xi(x) \]
is obtained. Since this equation holds for all $x\in \pi(E)$ and
$\pi_*(\nu)(\pi(E))\geq\nu(E)=1$,
the proof is completed.
\end{proof}

Let $\mu\in M_\alpha$.
As in the discussion following Definition \ref{KT},
let $\xi(x)=\theta$ be a constant function.
Then, there exists a Borel function $\eta\in C(X,\T)$ such that
\[ \xi(x)=\eta(x)-\eta\alpha^{-1}(x) \]
for $\mu$-almost every $x\in X$ if and only if
$e^{2\pi\sqrt{-1}\theta}$ is an eigenvalue of the unitary operator
$\pi_\mu(u_\alpha)\in L^2(X,\mu)$,
where $\pi_\mu$ is a representation of $C^*\xa$
corresponding to the invariant measure $\mu$.
Since $L^2(X,\mu)$ is separable, eigenvalues of a unitary operator
are at most countable.
Hence, by the lemma above,
we can obtain a lot of rigid homeomorphisms.
Moreover, it is known that eigenvalues of the unitary operator
$\pi_\mu(u_\alpha)$ need not be topological eigenvalues of $\xa$.
Therefore we can see that
there exists a minimal homeomorphism $\alpha\times R_\xi$
which is not rigid.
We will look at a concrete example in Example \ref{Furstenberg}.

There is another way to find a rigid homeomorphism.
Let $\xi\in C(X,\T)$ and
let $\tilde{\xi}\in C(X,\R)$ be its lift
(this is always possible because $X$ is the Cantor set).
Then $\mu \mapsto \mu(\tilde{\xi})$ gives an affine function
from the set of $\alpha$-invariant measures $M_\alpha$ to $\R$.
By the lemma above, if $n\mu(\tilde{\xi})\notin \mu(C(X,\Z))$
for each ergodic $\alpha$-invariant measure $\mu$ and $n\in\N$,
then $\alpha\times R_\xi$ is rigid.
\bigskip

Next, we will show that $A=C^*(X\times\T,\alpha\times R_\xi)$ can be
written as a crossed product of $C^*\xa$ by a certain action.
%
Define an automorphism $\iota(\xi)$ on $C^*\xa$ by
$\iota(\xi)(f)=f$ for all $f\in C(X)$ and
$\iota(\xi)(u_\alpha)=u_\alpha e^{2\pi\sqrt{-1}\xi(x)}$,
where $u_\alpha$ denotes the implementing unitary of $C^*\xa$.
This kind of automorphism was considered in \cite{M1}.
We remark that $\iota(\xi)$ is approximately inner,
because $\iota(\xi)_*$ is the identity on the $K$-groups
(or one can deduce it from Lemma \ref{perturb} or
\cite[Lemma 5.1]{M1}).
Let $\hat{\iota}(\xi)$ denote the dual action
on $C^*\xa\rtimes_{\iota(\xi)}\Z$.

\begin{prop}\label{ATrtimesZ}
There is an isomorphism $\pi$
from the crossed product $C^*$-algebra
$A=C^*(X\times\T,\alpha\times R_\xi)$ to $C^*\xa\rtimes_{\iota(\xi)}\Z$
such that $\pi(f)=f$ for all $f\in C(X)$ and
$\pi(g\circ(\id\times R_t))=\hat{\iota}(\xi)_t(\pi(g))$
for all $g\in C(X\times\T)$ and $t\in\T$.
\end{prop}
\begin{proof}
In order to avoid confusion, we have to use different symbols for
three implementing unitaries:
we denote the implementing unitary in $C^*\xa$ by $u_\alpha$ and
denote the unitary implementing $\iota(\xi)$ by $v$,
while $u\in A$ denotes the unitary implementing $\alpha\times R_\xi$.

Let $z\in C(X\times\T)$ be a unitary defined by
$z(x,t)=e^{2\pi\sqrt{-1}t}$.
Define $\pi(z)=v$ and $\pi(f)=f$
for all $f\in C(X)\subset C(X\times\T)$.
This is well-defined because $v$ and $f$ commute in $C^*\xa\rtimes\Z$.
Moreover, $\pi$ is an isomorphism from $C(X\times\T)$ onto its image.
We define $\pi(u)=u_\alpha$.
It is easy to check that
\[ \pi(u)\pi(f)\pi(u^*)=f\alpha^{-1}=\pi(f\alpha^{-1}) \]
for all $f\in C(X)$ and that
\[ \pi(u)\pi(z)\pi(u^*)=e^{-2\pi\sqrt{-1}\xi(\alpha^{-1}(x))}v
=\pi(e^{-2\pi\sqrt{-1}\xi(\alpha^{-1}(x))}z)=\pi(z\circ\alpha^{-1}). \]
Therefore $\pi$ is a homomorphism from $A$ to $C^*\xa\rtimes\Z$.
Clearly $\pi$ is surjective.
It is also straightforward to see
$\pi(g\circ(\id\times R_t)^{-1})=\hat{\iota}(\xi)_t(\pi(g))$
for all $g\in C(X\times\T)$ and $t\in\T$.

It remains to show that $\pi$ is an isomorphism.
Let $E$ be the conditional expectation from $A$ to $C(X\times\T)$.
It is well-known that $E$ is faithful.
We can define an action of $\T$ on $C^*\xa\times\Z$ by
\[ \gamma_t(f)=f, \ \gamma_t(u_\alpha)=e^{2\pi\sqrt{-1}t}u_\alpha, \
\text{ and } \ \gamma_t(v)=v \]
for $t\in\T$.
Let
\[ E_0(a)=\int_\T \gamma_t(a)\,dt \]
for $a\in C^*\xa\rtimes\Z$.
Then we have $\pi\circ E=E_0\circ\pi$.
The faithfulness of $E$ leads us to the conclusion.
\end{proof}

\begin{rem}\label{states}
By \cite[Corollary 5.7]{Pu1} or \cite[Theorem 5.5]{HPS}
there exist bijective correspondences
between the following spaces.
\begin{enumerate}
\item The state space $S(K^0\xa)$ of $K^0\xa$.
\item The tracial state space $T(C^*\xa)$ of $C^*\xa$.
\item The set $M_\alpha$ of
all $\alpha$-invariant probability measures.
\end{enumerate}
By Proposition \ref{ATrtimesZ} and Lemma \ref{whenrigid},
(2) and (3) of the above are also identified
with the following.
\begin{enumerate}
\item[(4)] The set of $\hat{\iota}(\xi)$-invariant traces
on $A=C^*(X\times\T,\alpha\times R_\xi)$.
\item[(5)] The set of probability measures on $X\times\T$
which are invariant under $\alpha\times R_\xi$ and
the rotation $\id\times R_t$.
\end{enumerate}
Suppose that $s\in S(K_0(A))$ is a state
on the ordered group $K_0(A)$.
Then there is $\nu\in M_\alpha$ such that
$s([f])=S_\nu([f])$ for all $[f]\in K^0\xa$,
where $S_\nu$ is a state on $K^0\xa$ coming from $\nu$
and $K^0\xa$ is viewed as a subgroup of $K_0(A)$.
The $\alpha$-invariant measure $\nu$ extends to
an $\alpha\times R_\xi$-invariant measure on $X\times\T$, and so
we can extend $S_\nu$ on $K_0(A)$
(different choices of $\alpha\times R_\xi$-invariant measures
do not concern the extension of $S_\nu$).
For $x\in K_0(A)$,
\[ S(K^0\xa)\ni S_\mu\mapsto S_\mu(x)\in\R \]
is an affine function on the state space $S(K^0\xa)$.
Since the image of $K^0\xa$ is dense in $\Aff(S(K^0\xa))$,
for any $\varepsilon>0$, there exists
$f_1,f_2\in C(X,\Z)$ such that
\[ S_\mu(x)-\varepsilon<S_\mu([f_1])<S_\mu(x)
<S_\mu([f_2])<S_\mu(x)+\varepsilon \]
for all $\mu\in M_\alpha$.
If $A$ is simple, then it follows from Theorem \ref{BSFC}
that $[f_1]<x<[f_2]$ in $K_0(A)$.
In particular, we have
\[ S_\nu(x)-\varepsilon<S_\nu([f_1])=s([f_1])<s(x)
<s([f_2])=S_\nu([f_2])<S_\nu(x)+\varepsilon. \]
Hence $s$ is equal to $S_\nu$ as a state on $K_0(A)$.
Consequently, the state space $S(K^0\xa)$ can be identified with
\begin{enumerate}
\item[(6)] The state space $S(K_0(A))$ of $A$
\end{enumerate}
when $A$ is simple.
\end{rem}


\begin{thm}
Let $\xa$ be a Cantor minimal system and
let $\xi\in C(X,\T)$.
Suppose that $\alpha\times R_\xi$ is a minimal homeomorphism
on $X\times\T$.
For the unital simple $C^*$-algebra
$A=C^*(X\times\T,\alpha\times R_\xi)$,
the following conditions are equivalent.
\begin{enumerate}
\item $\alpha\times R_\xi$ is rigid.
\item $A$ has real rank zero.
\item For every extremal tracial state $\tau$ on
$C^*\xa$ and every $n\in\N$,
$\iota(\xi)^n$ is not weakly inner
in the GNS representation $\pi_\tau$.
\end{enumerate}
If $C^*\xa$ has finitely many extremal traces,
then the conditions above are also equivalent to
\begin{enumerate}
\setcounter{enumi}{3}
\item $\iota(\xi)$ has the tracial cyclic Rohlin property.
\end{enumerate}
\end{thm}
\begin{proof}
(1)$\Leftrightarrow$(2) was shown in Corollary \ref{rr0}.

(1)$\Rightarrow$(3) follows \cite[Proposition 2.3]{K}.

(3)$\Rightarrow$(1).
Suppose that $\alpha\times R_\xi$ is not rigid.
By Lemma \ref{whenrigid},
there exist an ergodic measure $\mu\in M_\alpha$,
$n\in\N$ and a Borel function $\eta:X\rightarrow\T$ such that
\[ n\xi(x)=\eta(x)-\eta\alpha(x) \]
for $\mu$-almost every $x\in X$.
Define $h\in L^{\infty}(X,\mu)$
by $h(x)=e^{2\pi\sqrt{-1}\eta(x)}$
and let $V$ be the multiplication operator
by $h$ on $L^2(X,\mu)$.
Let $\tau$ be the extremal trace on $C^*\xa$
corresponding to $\mu$.
We can regard $\pi_\tau(C^*\xa)$ as a $C^*$-subalgebra of
$B(L^2(X,\mu))$.
Then, it is not hard to see that
$V$ commutes with $\pi_\tau(f)$ for all $f\in C(X)$ and
that
\[ V^*\pi_\tau(u_\alpha)V=
\pi_\tau(\iota(\xi)^n(u_\alpha)). \]
Namely $\iota(\xi)^n$ is weakly inner
in the GNS representation $\pi_\tau$.

(3)$\Rightarrow$(4).
From \cite{OP} we can see that
$\iota(\xi)$ has the tracial Rohlin property.
The conclusion follows from \cite[Theorem 3.4]{LO}.

(4)$\Rightarrow$(2).
It follows from \cite[Theorem 2.9]{LO} that $A$ has tracial rank zero.
The conclusion is immediate from \cite[Theorem 3.4]{L1}.
\end{proof}

In the theorem above, it was shown that
if $\alpha\times R_\xi$ is rigid and
$\alpha$ has only finitely many ergodic measures,
then $A$ has tracial rank zero.
In the next section we will prove that the hypothesis of
finitely many ergodic measures is actually not necessary.

\section{Tracial rank}

Throughout this section, let $\xa$ be a Cantor minimal system and
let $\xi:X\rightarrow\T$ be a continuous map.
We denote the crossed product $C^*$-algebra
$C^*(X\times\T,\alpha\times R_\xi)$ by $A$ and
its implementing unitary by $u$.
We would like to show that if $\alpha\times R_\xi$ is rigid
then $A$ has tracial rank zero.
The proof will be done by some improvement of Lemma \ref{almostTAF}.
Following the notation used there, we define
$A_x=C^*(C(X\times\T),uC_0((X\setminus\{x\})\times\T))$ for $x\in X$.

Define $z\in C(X\times\T)$ by $z(x,t)=e^{2\pi\sqrt{-1}t}$.
The key step of the proof is approximately unitary equivalence of
$z1_U$ and $z1_V$ in $A_x$,
where $U$ and $V$ are suitable clopen subsets of $X$
satisfying $[1_U]=[1_V]$.
When one uses the fact that $A_x$ has tracial rank zero,
the proof is just an application of \cite[Theorem 3.4]{Lnd},
in which actually more general result has been obtained.
But we would like to include an elementary proof which does not use
tracial rank zero for the reader's convenience.

The following lemma says that rigidity implies that
the values of a cocycle are uniformly distributed in $\T$.

\begin{lem}\label{ufmdtrb}
Suppose that $\alpha\times R_\xi$ is rigid.
For any irrational $s\in\T$ and any $\varepsilon>0$,
there exists $N\in\N$ such that the following is satisfied:
for any $n\geq N$ and $y\in X$ there is a permutation $\sigma$
on $\{1,2,\dots,n\}$ such that
\[ \lvert ks-\sum_{i=0}^{\sigma(k)-1}\xi(\alpha^i(y))\rvert
<\varepsilon \]
holds for all $k\in\{1,2,\dots,n\}$.
\end{lem}
\begin{proof}
For $n\in\N$, put
$\xi_n(x)=\sum_{i=0}^{n-1}\xi(\alpha^i(x))$.
By Lemma \ref{whenrigid}, we have
\[ \int_{X\times\T}f(t)\,d\nu=\int_\T f(t)\,dt \]
for every invariant measure $\nu$ of $\alpha\times R_\xi$ and
every $f\in C(\T)$.
Hence for any $f\in C(\T)$ and $\varepsilon>0$
there exists $N\in\N$ such that
\[ \left\lvert\frac{1}{n}\sum_{i=0}^{n-1}f(\xi_n(x))-
\int_\T f(t)\,dt\right\rvert<\varepsilon \]
for all $n\geq N$ and $x\in X$.
By a slight modification of \cite[Lemma 2]{KK},
we can get the conclusion.
We leave the details to the reader.
\end{proof}

In the following lemmas we need the idea of induced transformations.
Let $U$ be a clopen subset of $X$. Define $r:U\rightarrow\N$ by
\[ r(x)=\min\{n\in\N:\alpha^n(x)\in U\}. \]
Since $\alpha$ is minimal, $r$ is well-defined and continuous.
Put $\tilde{\alpha}(y)=\alpha^{r(y)}(y)$ for every $y\in U$.
Thus $\tilde{\alpha}$ is the first return map on $U$.
It is well-known that $(U,\tilde{\alpha})$ is a Cantor minimal system
and the associated crossed product $C^*(U,\tilde{\alpha})$ is
canonically identified with $1_UC^*\xa1_U$.
Define $\tilde{\xi}:U\rightarrow\T$ by
\[ \tilde{\xi}(y)=\sum_{i=0}^{r(y)-1}\xi(\alpha^i(y)) \]
for all $y\in U$.
Then $\tilde{\alpha}\times R_{\tilde{\xi}}$ is the first return map
of $\alpha\times R_\xi$ on $U\times\T$ and
the associated crossed product
$C^*(U\times\T,\tilde{\alpha}\times R_{\tilde{\xi}})$ is identified with
$1_{U\times\T}A1_{U\times\T}$.
Note that the unitary implementing $\tilde{\alpha}\times R_{\tilde{\xi}}$
is given by
\[ \sum_{n\in\N}u^n1_{U_n\times\T}, \]
where $U_n=r^{-1}(n)$ and the summation is actually finite.

In general there is a bijective correspondence
between invariant measures of the induced transformation and
those of the original one.
It follows that $\alpha\times R_\xi$ is rigid if and only if
$\tilde{\alpha}\times R_{\tilde{\xi}}$ is rigid.

For $x\in X$, let $k$ be the minimal natural number such that
$\alpha^{-k}(x)\in U$ and set $\tilde{x}=\alpha^{-k}(x)$.
Then it is not hard to see that
\[ 1_UA_x1_U=C^*(C(U\times\T),
\tilde{u}C_0((U\setminus\{\tilde{x}\})\times\T)). \]

\begin{lem}\label{approxunitary}
Let $\alpha\times R_\xi$ be a rigid homeomorphism
and let $x\in X$.
Suppose that $U$ is a nonempty clopen subset of $X$.
For any $s\in\T$ and any $\varepsilon>0$,
there exists a unitary $w\in 1_U(A_x\cap C^*\xa)1_U$ such that
\[ \lVert wzw^*-e^{2\pi\sqrt{-1}s}z1_U\rVert<\varepsilon. \]
\end{lem}
\begin{proof}
At first we consider the case $U=X$.
Clearly we may assume that $s$ is irrational.
By applying Lemma \ref{ufmdtrb} we can find $N\in\N$.
Let
\[ {\cal P}=\{X(v,k):v\in V,k=1,2,\dots,h(v)\} \]
be a Kakutani-Rohlin partition such that the roof set
$R({\cal P})$ contains $x$ and
$h(v)$ is greater than $N$ for every $v\in V$.
By dividing each tower if necessary,
we may assume that ${\cal P}$ is sufficiently finer so that
whenever $y_1,y_2\in X(v,k)$ we have
$\lvert\xi(y_1)-\xi(y_2)\rvert<\varepsilon/h(v)$.

For every $v\in V$, choose $y_v\in X(v,1)$ arbitrarily.
By Lemma \ref{ufmdtrb} there is a permutation $\sigma_v$
on $\{1,2,\dots,h(v)\}$ such that
\[ \lvert ks-\sum_{i=0}^{\sigma_v(k)-1}\xi(\alpha^i(y_v))\rvert
<\varepsilon \]
for all $k\in\{1,2,\dots,h(v)\}$.
Put
\[ w=\sum_{v\in V}\sum_{i=1}^{h(v)}
1_{X(v,\sigma_v(i))}u^{\sigma_v(i)-\sigma_v(i+1)}, \]
where $u$ is the implementing unitary of $C^*\xa$.
It is easily verified that $w$ is a unitary of $A_x\cap C^*\xa$.
Moreover we get the estimate
\[ \lVert wzw^*-e^{2\pi\sqrt{-1}s}z\rVert<4\pi\varepsilon. \]

Let us consider the general case.
We follow the notation used in the discussion before the lemma.
Applying the first part of the proof to
$\tilde{\alpha}\times\tilde{\xi}$ and $\tilde{x}$,
we obtain a unitary $w$ in
\[ C^*(C(U),\tilde{u}C_0(U\setminus\{\tilde{x}\}))=
1_U(A_x\cap C^*\xa)1_U \]
which satisfies the required inequality.
\end{proof}

\begin{lem}\label{approxunitaryII}
Let $\alpha\times R_\xi$ be a rigid homeomorphism
and let $x\in X$.
For any $\eta\in C(X,\T)$ and any $\varepsilon>0$,
there exists a unitary $w\in A_x\cap C^*\xa$ such that
\[ \lVert wzw^*-e^{2\pi\sqrt{-1}\eta}z\rVert<\varepsilon, \]
where $z\in C(X\times\T)$ is given by $z(x,t)=e^{2\pi\sqrt{-1}t}$.
\end{lem}
\begin{proof}
Let ${\cal P}$ be a partition of $X$ such that
whenever $y_1,y_2\in U\in{\cal P}$ we have
$\lvert\eta(y_1)-\eta(y_2)\rvert<\varepsilon$.
For every $U\in{\cal P}$, by Lemma \ref{approxunitary},
we obtain a unitary $w_U\in1_U(A_x\cap C^*\xa)1_U$ satisfying
\[ \lVert w_Uzw_U^*-
e^{2\pi\sqrt{-1}\eta(\tilde{x})}z1_U\rVert<\varepsilon. \]
Let $w$ be the product of all $w_U$'s.
Then $w$ is the desired unitary.
\end{proof}

\begin{lem}\label{approxunitaryIII}
Let $\alpha\times R_\xi$ be a rigid homeomorphism
and let $x\in X$.
Suppose that $U$ is a clopen neighborhood of $x$ and
$U,\alpha(U),\dots,\alpha^M(U)$ are mutually disjoint.
Put $p=1_U$ and $q=1_{\alpha^M(U)}$.
Then for any $\varepsilon>0$
there exists a partial isometry $w\in A_x\cap C^*\xa$
such that $w^*w=p$, $ww^*=q$ and
\[ \lVert wzw^*-zq\rVert<\varepsilon. \]
Moreover we have $u^{*i}wu^i\in A_x\cap C^*\xa$
for all $i=0,1,\dots,M-1$.
\end{lem}
\begin{proof}
There exists a partial isometry $v_1\in A_x\cap C^*\xa$
such that $v_1^*v_1=p$ and $v_1v_1^*=q$.
We have $v_1^*zv_1=e^{2\pi\sqrt{-1}\eta}zp$
for some continuous function $\eta$ defined on $U$.
We consider the induced transformation on $U$.
Let $\tilde{\alpha}$, $\tilde{\xi}$ and $\tilde{x}$ be
as in the discussion before Lemma \ref{approxunitary}.
Then Lemma \ref{approxunitaryII} applies to them and yields
a unitary $v_2\in p(A_x\cap C^*\xa)p$ satisfying
\[ \lVert v_2zpv_2^*-e^{2\pi\sqrt{-1}\eta}zp\rVert<\varepsilon. \]
Then $w=v_1v_2$ satisfies
\[ \lVert wzw^*-zq\rVert<\varepsilon. \]
Since $U,\alpha(U),\dots,\alpha^M(U)$ are mutually disjoint,
one can check that $w$ belongs to $A_{\alpha^i(x)}\cap C^*\xa$
for all $i=0,1,\dots,M-1$.
It follows that $u^{*i}wu^i\in A_x\cap C^*\xa$
for all $i=0,1,\dots,M-1$.
\end{proof}

\begin{lem}\label{Lad1}
Suppose that $\alpha\times R_\xi$ is rigid. Let $x\in X$.
For any $N\in\N$, $\varepsilon>0$ and
a finite subset ${\cal F}\subset C(X\times\T)$,
we can find a natural number $M>N$,
a clopen neighborhood $U$ of $x$ and a partial isometry $w\in A_x$
which satisfy the following.
\begin{enumerate}
\item $\alpha^{-N+1}(U),\alpha^{-N+2}(U),\dots,
U,\alpha(U),\dots,\alpha^M(U)$
are mutually disjoint, and
$\mu(U)<\varepsilon/M$ for all $\alpha$-invariant measure $\mu$.
\item $w^*w=1_U$ and $ww^*=1_{\alpha^M(U)}$.
\item $u^{*i}wu^i\in A_x$ for all $i=0,1,\dots,M-1$.
\item $\lVert wf-fw\rVert<\varepsilon$ for all $f\in{\cal F}$.
\end{enumerate}
\end{lem}
\begin{proof}
Without loss of generality, we may assume
${\cal F}=\{f_1,f_2,\dots,f_k,z\}$,
where $f_i$ belongs to $C(X)\subset C(X\times\T)$.
There exists a clopen neighborhood $O$ of $x$ such that
\[ \lvert f_i(x)-f_i(y)\rvert<\varepsilon/2 \]
for all $y\in O$ and $i=1,2,\dots,k$.
Since $\alpha$ is minimal, we can find $M>N$ such that
$\alpha^M(x)\in O$.
Let $U$ be a clopen neighborhood of $x$ such that
the condition (1) is satisfied and $U\cup\alpha^M(U)\subset O$.
Now Lemma \ref{approxunitaryIII} applies and yields
a partial isometry $w$.
It is clear that $w$ is the desired one.
\end{proof}

\begin{thm}\label{Tad1}
Suppose that $\alpha\times R_\xi$ is minimal.
Then the following are equivalent.
\begin{enumerate}
\item $\alpha\times R_\xi$ is rigid.
\item $A=C(X\times\T,\alpha\times R_\xi)$ has real rank zero.
\item $A=C(X\times\T,\alpha\times R_\xi)$ has tracial rank zero.
\item $A=C(X\times\T,\alpha\times R_\xi)$ is a unital simple AT-algebra
with real rank zero.
\end{enumerate}
\end{thm}
\begin{proof}
It has been proved in Corollary \ref{rr0} that (1) and (2) are equivalent.

(3)$\Leftrightarrow$(4) follows the classification theorem
of \cite{L3}, since $A$ has torsion free $K$-theory
(see Lemma \ref{Kofop}).

(4)$\Rightarrow$(2) is obvious.

(1)$\Rightarrow$(3).
We will show the following:
For any $\varepsilon>0,$
any finite subset ${\cal F}\subset C(X\times\T)$
and any nonzero positive element $c\in A,$
there exists a projection $e\in A_x$
such that the following conditions hold.
\begin{itemize}
\item $\lVert ae-ea\rVert<\varepsilon$
for all $a\in {\cal F}\cup\{u\}$.
\item For any $a\in {\cal F}\cup\{u\}$, there exists $b\in eA_xe$
such that $\lVert eae-b\rVert<\varepsilon$.
\item $1-e$ is equivalent to a projection in $\overline{cAc}.$
\end{itemize}
It follows from Proposition \ref{Putcutdown}
that $eA_xe$ is a unital simple AT-algebra with real rank zero.
Therefore it has tracial rank zero
(for example see \cite[Theorem 4.3.5]{Lnb}).
Thus, if the above is proved, by \cite[Theorem 4.8]{HLY},
$A$ has tracial rank zero.
In Section 3 it was proved that $A$ has real rank zero, stable
rank one and has weakly unperforated $K_0(A),$
and so it suffices to show the following:
For any $\varepsilon>0$ and
a finite subset ${\cal F}\subset C(X\times\T)$,
there exists a projection $e\in A_x$ such that
the following conditions hold.
\begin{itemize}
\item $\lVert ae-ea\rVert<\varepsilon$
for all $a\in {\cal F}\cup\{u\}$.
\item For any $a\in {\cal F}\cup\{u\}$, there exists $b\in eA_xe$
such that $\lVert eae-b\rVert<\varepsilon$.
\item $\tau(1-e)<\varepsilon$ for all $\tau\in T(A)$.
\end{itemize}

We may assume ${\cal F}^*={\cal F}$.
Choose $N\in\N$ so that $2\pi/N$ is less than $\varepsilon$.
Applying Lemma \ref{Lad1} to $N$, $\varepsilon/2$ and a finite subset
\[ {\cal G}=\bigcup_{i=0}^{N-1}u^i{\cal F}u^{i*}, \]
we obtain $M>N$, a clopen neighborhood $U$ of $x$
and a partial isometry $w\in A_x$.
Put $p=1_U$ and $q=1_{\alpha^M(U)}$.
For $t\in[0,\pi]$ we define
\[ P(t)=p\cos^2 t+w\sin t\cos t+w^*\sin t\cos t+q\sin^2 t. \]
Then $P(t)$ is a continuous path of projections with $P(0)=p$ and
$P(\pi)=q$.
By the choice of $w$ we obtain the estimate
\[ \lVert u^{i*}P(t)u^if-fu^{i*}P(t)u^i\rVert<\varepsilon \]
for all $t\in[0,\pi]$, $i=0,1,\dots,N-1$ and $f\in {\cal F}$.
We define a projection $e$ by
\[ e=1-\left(\sum_{i=0}^{M-N}u^ipu^{i*}+
\sum_{i=1}^{N-1}u^{i*}P(i\pi/N)u^i\right). \]
The partial isometry $w$ satisfies $u^{i*}wu^i\in A_x$
for all $i=1,2,\dots,N-1$, and so $e$ is a projection of $A_x$.
Evidently we have
\[ \lVert fe-ef\rVert<\varepsilon \]
for all $f\in{\cal F}$.
Since
\[ \lVert P(i\pi/N)-P((i-1)\pi/N)\rVert<
\frac{2\pi}{N}<\varepsilon, \]
it is not hard to see
\[ \lVert ue-eu\rVert<\varepsilon. \]
It is clear that $efe$ belongs to $A_x$ for all $f\in C(X\times\T)$.
It follows from $eue=eu(1-p)e$ that $eue$ also belongs to $A_x$.
We can easily verify
\[ \tau(1-e)<M\tau(p)<\frac{\varepsilon}{2} \]
for all $\tau\in T(A)$.
\end{proof}

\section{The generalized Rieffel projection}

By Lemma \ref{Kofop}, the $K_i$-group ($i=1,2$) of the crossed product
$C^*$-algebra arising from an orientation preserving
homeomorphism $\alpha\times\phi$ is isomorphic to the direct sum
of $\Z$ and $K^0\xa$.
Needless to say, the equivalence class of the implementing unitary is
the generator of $\Z$ in the $K_1$-group.
This section is devoted to specify a projection of
$C^*(X\times\T,\alpha\times\phi)$ which gives a representative
of the generator of the $\Z$-summand of the $K_0$-group.

At first, we consider the case that a cocycle takes its values
in the rotation group of the circle.
Let $\xa$ be a Cantor minimal system and
let $\xi:X\rightarrow\T$ be a continuous map.
We denote $C^*(X\times\T,\alpha\times R_\xi)$ by $A$ for short and
denote the implementing unitary by $u$.
We will identify $K^0\xa$ as a subgroup of $K_0(A)$.

Let $\tilde{\xi}\in C(X,\R)$ be an arbitrary lift of
$\xi\in C(X,\T)$. Then $M_\alpha\ni\mu\mapsto
\mu(\tilde{\xi})\in\R$
gives an affine function
on the set of invariant probability measures $M_\alpha$.
The other lifts of $\xi$ are of the form
$\tilde{\xi}+f$ with $f\in C(X,\Z)$, and so
this affine function is uniquely determined up to
the natural image of $K^0\xa$ in $\Aff(M_\alpha)$.
Suppose $[\xi]=0$ in $K^0_\T\xa$.
Then there exists $\eta\in C(X,\T)$ such that
$\xi=\eta-\eta\alpha^{-1}$.
When $\tilde{\eta}\in C(X,\R)$ is a lift of $\eta$,
$\tilde{\xi}=\tilde{\eta}-\tilde{\eta}\alpha^{-1}$
is a lift of $\xi$ and $\mu(\tilde{\xi})=0$
for all $\mu\in M_\alpha$.
Therefore we obtain a homomorphism
\[ \Delta:K^0_\T\xa\ni[\xi] \ \mapsto \
\Delta([\xi])\in\Aff(M_\alpha)/D(K^0\xa). \]
This homomorphism $\Delta$ is known to be surjective
(\cite[Lemma 6.2]{M1}).

\begin{lem}\label{perturb}
Let $\xa$ be a Cantor minimal system and
let $\xi:X\rightarrow\T$ be a continuous map.
For any $\varepsilon>0$, there exists
$\eta\in C(X,\T)$ such that
$\lvert(\xi+\eta-\eta\alpha)(x)\rvert<\varepsilon$
for all $x\in X$.
\end{lem}
\begin{proof}
Let
\[ {\cal P}=\{X(v,k):v\in V,k=1,2,\dots,h(v)\} \]
be a Kakutani-Rohlin partition of $\xa$ such that
$h(v)>\varepsilon^{-1}$ for every $v\in V$.
We denote the roof set $\bigcup_{v\in V}X(v,h(v))$
by $R({\cal P})$.
Put
\[ \kappa(x)=\sum_{k=1}^{h(v)}\xi(\alpha^{k-1}(x)) \]
for all $x\in X(v,1)$ and $v\in V$.
Since $X$ is totally disconnected,
there exists a real valued continuous function
$\tilde{\kappa}$ on $\alpha(R({\cal P}))$ such that
\[ \tilde{\kappa}(x)+\Z=\kappa(x) \ \text{ and } \
-1<\tilde{\kappa}(x)<1 \]
for all $x\in\alpha(R({\cal P}))=\bigcup_{v\in V}X(v,1)$.
Define $\eta\in C(X,\T)$ by
$\eta(x)=0$ for all $x\in\alpha(R({\cal P}))$ and
\[ \eta(\alpha^k(x))=\sum_{i=1}^{k}\xi(\alpha^{i-1}(x))-
\frac{k}{h(v)}\tilde{\kappa}(x)+\Z \]
for $x\in X(v,1)$, $v\in V$ and $k=1,2,\dots,h(v)-1$.
It is not hard to see that $\eta$ is the desired function.
\end{proof}

In a similar fashion to the lemma above, we can show the following,
which will be used later.

\begin{lem}\label{perturbIII}
Let $\xa$ be a Cantor minimal system and
let $\xi:X\rightarrow\T$ and $c:X\rightarrow\Z_2$ be continuous maps.
For any $\varepsilon>0$, there exists
$\eta\in C(X,\T)$ such that
\[ \lvert\xi(x)+\eta(x)-(-1)^{c(x)}\eta\alpha(x)\rvert<\varepsilon \]
for all $x\in X$.
\end{lem}

\begin{df}
Let $\xa$ be a Cantor minimal system and
let $\xi:X\rightarrow\T$ be a continuous map.
Define
\[ H(\alpha,\xi)=\{\eta\in C(X,\T) \ : \
(\xi+\eta-\eta\alpha)(x)\in(1/10,9/10)
\text{ for all }x\in X\}. \]
By Lemma \ref{perturb}, $H(\alpha,\xi)$ is not empty.
\end{df}

Suppose that $\eta$ belongs to $H(\alpha,\xi)$.
We define a projection $e(\alpha,\xi,\eta)$
in $A$ as follows.
Define a real valued continuous function $g_{\eta}$ on
$X\times\T$ by
\[ g_\eta(x,t)=\begin{cases}
\sqrt{10(t-\eta(x))(1-10(t-\eta(x)))} &
t\in [\eta(x),\eta(x)+1/10] \\
0 & \text{otherwise.} \end{cases} \]
Put $\eta'=(\xi+\eta)\circ\alpha^{-1}$ and
define a real valued continuous function $f(\alpha,\xi,\eta)$
on $X\times\T$ by
\[ f(\alpha,\xi,\eta)(x,t)=\begin{cases}
10(t-\eta(x)) & t\in[\eta(x),\eta(x)+1/10] \\
1 & t\in[\eta(x)+1/10,\eta'(x)] \\
1-10(t-\eta'(x)) & t\in[\eta'(x),\eta'(x)+1/10] \\
0 & \text{otherwise.} \end{cases} \]
Then it is easy to check that
$e(\alpha,\xi,\eta)=g_\eta u^*+f(\alpha,\xi,\eta)+ug_\eta\in A$
is a self-adjoint projection.
We call $e(\alpha,\xi,\eta)$
the generalized Rieffel projection.
Let $\tilde{\xi}\in C(X,\R)$ be a lift of
$\xi+\eta-\eta\alpha$ such that $1/10<\tilde{\xi}<9/10$.
Then, for every $\mu\in M_\alpha$,
\[ \tau_\mu(e(\alpha,\xi,\eta))=
\tau_\mu(f(\alpha,\xi,\eta))=\mu(\tilde{\xi}), \]
where $\tau_\mu$ is the tracial state on $A$ corresponding to $\mu$.
Hence the affine function $\mu\mapsto\tau_\mu(e(\alpha,\xi,\eta))$
is a representative of
$\Delta([\xi])\in\Aff(M_\alpha)/D(K^0\xa)$.

\begin{prop}\label{RieffelI}
In the situation above,
let $e=e(\alpha,\xi,\eta)\in A$ be the generalized Rieffel
projection.
Then $K^0\xa$ and $[e]$ generate $K_0(A)$.
\end{prop}
\begin{proof}
Let $v$ the unilateral shift on $\ell^2(\N)$ and
let ${\cal T}$ be the Toeplitz algebra generated by $v$.
Put $q=1-vv^*$.
In the $C^*$-algebra $A\otimes{\cal T}$, we consider
the $C^*$-subalgebra $B$ generated
by $C(X\times\T)\otimes 1$ and $u\otimes v^*$.
There is a surjective homomorphism $\pi$ from $B$
to $A$ sending $u\otimes v^*$ to $u$.
The kernel of $\pi$ is $C(X\times\T)\otimes{\cal K}$,
where ${\cal K}$ is the algebra of compact operators.
Put $a=g_\eta u^*\otimes v+f(\alpha,\xi,\eta)\otimes 1
+ug_\eta\otimes v^*$.
Then $\pi(a)=e$.
Hence it suffices to show that $e^{2\pi\sqrt{-1}a}$ is
a generator of $K_1(\C1_X\otimes C(\T))\cong\Z$.
Define a continuous function $h(t)$ by
\[ h(t)=(1-8t(1-t))+\sqrt{-1}t(t-1)(t-1/2). \]
Since $e^{2\pi\sqrt{-1}t}$ is homotopic to $h(t)$
in the set of complex-valued invertible functions,
it suffices to know the $K_1$-class
of the invertible element $h(a)$.
By
\[ a-a^2=g_\eta^2\otimes q \]
and
\[ (a^2-a)(a-1/2)=-(g_\eta^2(f(\alpha,\xi,\eta)-1/2))\otimes q, \]
it follows that $h(a)$ is homotopic to
$1_X\otimes z\otimes q$ in $GL(C(X\times\T))\otimes q$,
where $z(t)=e^{2\pi\sqrt{-1}t}$ is the generator of $K_1(C(\T))$.
\end{proof}
Let $K_0(A)\cong\Z\oplus K^0\xa$ be the isomorphism described
in the proposition above.
If $\alpha\times R_\xi$ is minimal,
then it follows from Theorem \ref{BSFC} that
\[ K_0(A)^+\cong\{(n,[f]):\mu(n\tilde{\xi}+f)>0 \ \text{ for all }
\ \mu\in M_\alpha\}\cup\{0\}. \]
See also Remark \ref{states}.
\bigskip

We now turn to the general case.
Let $\xa$ be a Cantor minimal system and
let $\phi:X\rightarrow\Homeo^+(\T)$ be a continuous map.
We write the crossed product $C^*$-algebra arising from
$(X\times\T,\alpha\times\phi)$ by $A$ and
the implementing unitary by $u$.
In order to define the Rieffel projection in $A$,
we need some preparations.
For $\phi\in\Homeo^+(\T)$, let $r(\phi)\in\T$ denote the rotation number.
The reader may refer to \cite[Chapter 11]{KH}
for the definition and some elementary properties of $r(\phi)$.

\begin{lem}\label{prepaA}
Let $\phi\in\Homeo^+(\T)$.
The map $\T\ni t\mapsto r(R_t\phi)$ is a continuous surjection
from $\T$ to $\T$ of degree one.
\end{lem}
\begin{proof}
By the definition of the rotation number, we see that
the mapping $t\mapsto r(R_t\phi)$ is nondecreasing
as a real valued function.
It is clear that $r(R_t\phi)=0$ if and only if $t$ belongs to
\[ I=\{s-\phi(s):s\in\T\}. \]
Since $\phi$ is an orientation preserving homeomorphism,
$I$ is not the whole circle.
Thus, $I$ is a closed interval of the circle.
It follows that the map is a surjection of degree one.
\end{proof}

\begin{lem}\label{prepaB}
Let $X$ be the Cantor set and let $I\subset\T$ be an open subset.
When $\Phi:X\times\T\rightarrow\T$ is a continuous map and
$\Phi(x,\cdot)$ is surjective for every $x\in X$,
there exists a continuous map $\xi:X\rightarrow\T$ such that
$\Phi(x,\xi(x))\in I$ for all $x\in X$.
\end{lem}
\begin{proof}
By assumption, for each $x\in X$, there exists $t_x\in\T$ such that
$\Phi(x,t_x)\in I$.
The continuity of $\Phi$ implies that
there exists a clopen neighborhood $U_x$ of $x$ such that
$\Phi(U_x,t_x)\subset I$.
Since $X$ is compact, it is covered by finitely many $U_x$'s.
We can find a locally constant function $\xi\in C(X,\T)$
satisfying the required property.
\end{proof}

The following lemma corresponds to Lemma \ref{perturb}.

\begin{lem}\label{perturbII}
Let $\xa$ be a Cantor minimal system and
let $\phi:X\rightarrow\Homeo^+(\T)$ be a continuous map.
Suppose that an open subset $I\subset\T$ is given.
Then there exists $\eta\in C(X,\T)$ such that
\[ r(R_{\eta\alpha(x)}\circ\phi_x\circ R_{\eta(x)}^{-1})\in I \]
for all $x\in X$.
\end{lem}
\begin{proof}
Define $\Phi(x,t)=r(R_t\phi_x)$.
It is obvious that $\Phi$ is continuous.
By Lemma \ref{prepaA}, for each $x\in X$, $\Phi(x,\cdot)$ is
a continuous surjection.
It follows from Lemma \ref{prepaB} that
there exists a continuous map $\xi:X\rightarrow\T$ such that
$\Phi(x,\xi(x))\in I$ for all $x\in X$.
Moreover, there exists $\varepsilon>0$ such that, for all $x\in X$,
if $\lvert t-\xi(x)\rvert<\varepsilon$ then $\Phi(x,t)\in I$.
From Lemma \ref{perturb}, we can find $\eta\in C(X,\T)$ such that
\[ \lvert(\xi+\eta-\eta\alpha)(x)\rvert<\varepsilon \]
for all $x\in X$. Then we get
\begin{align*}
& r(R_{\eta\alpha(x)}\circ\phi_x\circ R_{\eta(x)}^{-1}) \\
&=r(R_{\eta(x)}\circ R_{\eta\alpha(x)-\eta(x)}\circ
\phi_x\circ R_{-\eta(x)}) \\
&=r(R_{\eta\alpha(x)-\eta(x)}\circ\phi_x) \\
&=\Phi(x,\eta\alpha(x)-\eta(x))\in I,
\end{align*}
thereby completing the proof.
\end{proof}

By the lemma above, without loss of generality,
we may always assume that $r(\phi_x)$ is not zero for all $x\in X$.
Put
\[ c=\inf\{\lvert\phi_x(t)-t\rvert,
\lvert\phi_x^{-1}(t)-t\rvert:(x,t)\in X\times\T\}. \]
Since $\phi_x$ has no fixed points, $c$ is a positive real number.
Take $s\in\T$.
We define a function $f(\alpha,\phi,s)\in C(X\times\T)$ by
\[ f(\alpha,\phi,s)(x,t)=\begin{cases}
c^{-1}(t-s) & t\in[s,s+c] \\
1 & t\in[s+c,\phi_{\alpha^{-1}(x)}(s)] \\
1-c^{-1}(\phi_{\alpha^{-1}(x)}^{-1}(t)-s) &
t\in[\phi_{\alpha^{-1}(x)}(s),\phi_{\alpha^{-1}(x)}(s+c)] \\
0 & \text{otherwise.} \end{cases} \]
By the choice of $c$, $f$ is well-defined.
Define a function $g_s\in C(X\times\T)$ by
\[ g_s(x,t)=\begin{cases}
\sqrt{c^{-1}(t-s)(1-c^{-1}(t-s))} & t\in[s,s+c] \\
0 & \text{otherwise.} \end{cases} \]
Then one checks that $e(\alpha,\phi,s)=g_su^*+f(\alpha,\phi,s)+ug_s$ is
a well-defined projection of $A$.
Let us call it the generalized Rieffel projection for $A$.
In exactly the same way as Proposition \ref{RieffelI},
we can show the following.

\begin{prop}\label{RieffelII}
In the above setting, $K_0(A)$ is generated by $K^0\xa$ and
$[e(\alpha,\phi,s)]$.
Furthermore, $\T\ni s\mapsto e(\alpha,\phi,s)\in A$ is
a continuous path of projections in $A$.
\end{prop}
In the definition of $e(\alpha,\phi,s)$, we can replace $g_s$
by $zg_s$, where $z$ is a complex number with $\lvert z\rvert=1$.
But, this choice does not matter to the homotopy equivalence class
of the projection.

\section{Approximate $K$-conjugacy}

Let us begin with recalling the definition of
weakly approximate conjugacy.
\begin{df}[{\cite[Definition 3.1]{LM}}]
Let $\xa$ and $\yb$ be dynamical systems on compact metrizable
spaces $X$ and $Y$.
We say that $\xa$ and $\yb$ are weakly approximately conjugate,
if there exist homeomorphisms $\sigma_n:X\rightarrow Y$ and
$\tau_n:Y\rightarrow X$ such that
$\sigma_n\alpha\sigma_n^{-1}$ converges to $\beta$ in $\Homeo(Y)$
and $\tau_n\beta\tau_n^{-1}$ converges to $\alpha$ in $\Homeo(X)$.
\end{df}
In \cite[Theorem 4.13]{LM}, it was shown that
two Cantor minimal systems are weakly approximately conjugate
if and only if they have the same periodic spectrum
(see also \cite[Theorem 3.1]{M3}).
Similar results were shown in \cite{M3} for dynamical systems
on the product of the Cantor set and the circle.

As one sees that in the above definition, there is no consistency
among $\sigma_n$ or $\tau_n.$ It is clear (see \cite{LM}) that
the relation can be made stronger if one requires some consistency
among $\sigma_n$ as well as $\tau_n.$ We hope such stronger
version of approximate conjugacy is more reasonable replacement
of (flip) conjugacy.

Suppose $\sigma_n\alpha\sigma_n^{-1}\rightarrow\beta$ in $\Homeo(Y)$.
In \cite[Proposition 3.2]{LM}, it was shown that
there exists an asymptotic morphism $\{\psi_n\}:B\rightarrow A$
such that
\[ \lim_{n\rightarrow\infty}\lVert\psi_n(f)-f\circ\sigma_n\rVert=0 \]
for all $f\in C(Y)$ and
\[ \lim_{n\rightarrow\infty}\psi_n(u_\beta)=u_\alpha, \]
where $u_\alpha$ and $u_\beta$ denote the implementing unitaries
in $C^*\xa$ and $C^*\yb$.
This observation, however, is far from the existence theorem
in classification theory,
which requires that $\{\psi_n\}$ carries an isomorphism of $K$-groups
(see \cite[Theorem 4.3]{L3} for instance).
As pointed out in \cite{LM}, we have to impose conditions
on the conjugating maps $\sigma_n$ so that the associated asymptotic
morphism has a nice property.
Taking account of this, we make the following definitions.
By an order and unit preserving homomorphism
$\rho:K_*(B)\rightarrow K_*(A)$,
we mean a pair of homomorphisms $\rho_i:K_i(A)\rightarrow K_i(B)$
($i=0,1$) such that $\rho_0([1_A])=[1_B]$ and
$\rho_0(K_0(A)^+)\subset K_0(B)^+$.

\begin{df}
Let $\xa$ and $\yb$ be dynamical systems on compact metrizable
spaces $X$ and $Y$.
Suppose that a sequence of homeomorphisms $\sigma_n:X\rightarrow Y$
satisfies $\sigma_n\alpha\sigma_n^{-1}\rightarrow\beta$ in $\Homeo(Y)$.
Let $\{\psi_n\}$ be the asymptotic morphism arising from $\sigma_n$.
We say that the sequence $\{\sigma_n\}$ induces
an order and unit preserving  homomorphism
$\rho:K_*(C^*\yb)\rightarrow K_*(C^*\xa)$ between $K$-groups,
if for every projection $p\in M_\infty(C^*\yb)$ and
every unitary $u\in M_\infty(C^*\yb)$,
there exists $N\in\N$ such that
\[ [\psi_n(p)]=\rho([p])\in K_0(C^*\xa) \ \text{ and } \
[\psi_n(u)]=\rho([u])\in K_1(C^*\xa) \]
for every $n\geq N$.
\end{df}

\begin{df}[{\cite[Definition 5.3]{Lnd}}]\label{appKconju}
Let $\xa$ and $\yb$ be dynamical systems on compact metrizable
spaces $X$ and $Y$.
We say that $\xa$ and $\yb$ are approximately $K$-conjugate,
if there exist homeomorphisms $\sigma_n:X\rightarrow Y$,
$\tau_n:Y\rightarrow X$ and
an isomorphism $\rho:K_*(C^*\yb)\rightarrow K_*(C^*\xa)$
between $K$-groups such that
\[ \sigma_n\alpha\sigma_n^{-1}\rightarrow\beta, \ \
\tau_n\beta\tau_n^{-1}\rightarrow\alpha \]
and the associated asymptotic morphisms $\{\psi_n\}:B\rightarrow A$
and $\{\phi_n\}:A\rightarrow B$ induce
the isomorphisms $\rho$ and $\rho^{-1}$.

We say that $\xa$ and $\yb$ are approximately flip $K$-conjugate,
if $\xa$ is approximately $K$-conjugate to
either of $\yb$ and $(Y,\beta^{-1})$.
\end{df}

J. Tomiyama \cite{Tm} proved that two topological transitive
systems $\xa$ and $\yb$ are flip conjugate if and only if
there is an isomorphism $\phi:C^*\xa\rightarrow C^*\yb$
such that $\phi\circ j_{\alpha}=j_{\beta}\circ \chi$
for some isomorphism $\chi:C(X)\to C(Y).$

\begin{df}[{\cite[Definition 3.8]{Lnd}}]\label{4dcapp}
Let $\xa$ and $\yb$ be two topological transitive systems.
We say that $\xa$ and $\yb$ are
$C^*$-strongly approximately flip conjugate
if there exists a sequence of isomorphisms
$\phi_n:C^*\xa\to C^\yb$ and
a sequence of isomorphisms $\chi_n:C(X)\to C(Y)$ such that
$[\phi_n]=[\phi_1]$ in $KL(C^*\xa,C^*\yb)$ for all $n\in\N$ and
$$
\lim_{n\to\infty}\lVert\phi_n\circ j_{\alpha}(f)-j_{\beta}\circ
\chi_n(f)\rVert=0
$$
for all $f\in C(X).$
\end{df}
\bigskip

Let $X$ be the Cantor set and
let $\xi$ be a continuous function from $X$ to $\T$.
In this section, we would like to discuss
approximate $K$-conjugacy for $(X\times\T,\alpha\times R_\xi)$.

As in the previous section, let $A$ denote the crossed product
$C^*$-algebra $C^*(X\times\T,\alpha\times R_\xi)$.
By Lemma \ref{perturb}, we can find a continuous function
$\zeta:X\rightarrow\T$ such that $[\zeta]=[\xi]$ in $K^0_\T\xa$ and
$\zeta(x)\in(7/15,8/15)$ for every $x\in X$.
By Lemma \ref{conju}, $\alpha\times R_\xi$ is conjugate to
$\alpha\times R_\zeta$.
Therefore we may assume $\xi(x)\in(7/15,8/15)$ for all $x\in X$
without loss of generality.

Suppose that $\eta$ belongs to $H(\alpha,\xi)$.
Let $\tilde{\eta}\in C(X,\R)$ be a lift of $\eta$.
Since $(\xi+\eta-\eta\alpha)(x)\in(1/10,9/10)$ for every $x\in X$,
we have $(\eta-\eta\alpha)(x)\notin [13/30,17/30]$.
Hence there exists a unique $f\in C(X,\Z)$ such that
\[ \lVert f-(\tilde{\eta}-\tilde{\eta}\alpha)\rVert_\infty
<\frac{1}{2}. \]
It is easy to see that $[f]\in K^0\xa$ is independent of
the choice of $\tilde{\eta}$.
Let us denote $[f]$ by $B_\alpha(\eta)$.

\begin{lem}\label{homotopic}
In the above setting,
suppose that $\eta$ and $\eta'$ are homotopic in $H(\alpha,\xi)$.
\begin{enumerate}
\item $e(\alpha,\xi,\eta)$ is homotopic to $e(\alpha,\xi,\eta')$
in the set of projections of $A$.
\item $B_\alpha(\eta)=B_\alpha(\eta')$.
\end{enumerate}
\end{lem}
\begin{proof}
(1) Suppose $[0,1]\ni t\mapsto \kappa_t\in H(\alpha,\xi)$ is a homotopy
from $\eta$ to $\eta'$.
Then the generalized Rieffel projection $e(\alpha,\xi,\kappa_t)$ is
well-defined and $t\mapsto e(\alpha,\xi,\kappa_t)$ gives
a continuous path of projections in $A$
from $e(\alpha,\xi,\eta)$ to $e(\alpha,\xi,\eta')$.

(2) Let $\kappa_t$ be a homotopy as in (1).
There exists a continuous map $\tilde{\kappa}$
from $X\times[0,1]$ to $\R$ such that
\[ \tilde{\kappa}_t(x)+\Z=\kappa_t(x) \]
for all $x\in X$ and $t\in[0,1]$. Since
\[ (\tilde{\kappa}_t-\tilde{\kappa}_t\alpha)(x)\neq\frac{1}{2}, \]
the integer nearest to
$(\tilde{\kappa}_t-\tilde{\kappa}_t\alpha)(x)$ does not vary
as $t$ varies.
Hence we get $B_\alpha(\eta)=B_\alpha(\eta')$.
\end{proof}

Let $\eta\in H(\alpha,\xi)$. The unitary
$v_\eta(x)=e^{2\pi\sqrt{-1}\eta(x)}$ of $C(X)$ satisfies
$\lVert u_\alpha v_\eta u_\alpha^*v_\eta^*-1\rVert_\infty<2$,
where $u_\alpha$ denotes the implementing unitary of $C^*\xa$,
because of $(\eta-\eta\alpha)(x)\neq1/2$.
Thus, $u_\alpha$ and $v_\eta$ are almost commuting unitaries in a sense.
When $u,v\in M_n(\C)$ are unitaries satisfying $\lVert uv-vu\rVert<2$,
on account of $\det(uvu^*v^*)=1$, we have
\[ \frac{1}{2\pi\sqrt{-1}}\Tr(\log (uvu^*v^*))
\in\Z\cong K_0(M_n(\C)), \]
where $\Tr$ is the standard trace on $M_n(\C)$ and
$\log$ is the logarithm with values in $\{z:\Im(z)\in(-\pi,\pi)\}$.
The Bott element for pairs of almost commuting unitaries
in a unital $C^*$-algebra is a generalization of this (see \cite{EL}).
More precisely, if $u$ and $v$ are unitaries in a unital $C^*$-algebra
and $\lVert uv-vu\rVert\approx0$, then a projection $B(u,v)$ and
an element of the $K_0$-group are obtained.
Our $B_\alpha(\eta)$ is just this $K_0$-class
for $u_\alpha$ and $v_\eta$.

\begin{lem}\label{Bott}
Suppose $(\xi+\eta-\eta\alpha)(x)\in (1/3,2/3)$ for all $x\in X$.
Then we have
\[ [e(\alpha,\xi,\eta)]=[e(\alpha,\xi,0)]-B_\alpha(\eta) \]
in $K_0(A)$.
\end{lem}
\begin{proof}
Put
\[ c=\inf\{\lvert(\xi+\eta-\eta\alpha)(x)\rvert:x\in X\}
-\frac{1}{3}. \]
Then $c$ is positive.
Choose a sufficiently finer Kakutani-Rohlin partition
\[ {\cal P}=\{X(v,k):v\in V,1\leq k\leq h(v)\} \]
for $\xa$.
Let $R({\cal P})$ be the roof set of ${\cal P}$.
We may assume that $h(v)$ is large and
\[ \sup\{\lvert\eta(x)-\eta(y)\rvert:x,y\in U\}<\frac{c}{2} \]
for every $U\in \widetilde{\cal P}$,
where
\[ \widetilde{\cal P}=\{X(v,k):v\in V,1\leq k<h(v)\}
\cup\{R({\cal P})\}. \]
Take $x_U\in U$ for all $U\in \widetilde{\cal P}$ and
define $\eta'\in C(X,\Z)$ by $\eta'(x)=\eta(x_U)$ for $x\in U$.
It is not hard to see that $\eta$ and $\eta'$ are homotopic
in $H(\alpha,\xi)$.
By Lemma \ref{homotopic}, $[e(\alpha,\xi,\eta)]=[e(\alpha,\xi,\eta')]$
and $B_\alpha(\eta)=B_\alpha(\eta')$.
Hence, by replacing $\eta$ by $\eta'$, we may assume that
$\eta$ is constant on each clopen set belonging to $\widetilde{\cal P}$.
Furthermore, by adding a constant function, we may assume that
$\eta(x)=0$ for all $x\in R({\cal P})$.

Note that $\lvert(\eta-\eta\alpha)(x)\rvert$ is less than
$2/3-7/15=1/5$ for all $x\in X$.
There is a unique lift $\tilde{\eta}\in C(X,\R)$ of $\eta$
such that $\tilde{\eta}(x)=0$ for all $x\in R({\cal P})$ and
$\lvert(\tilde{\eta}-\tilde{\eta}\alpha)(y)\rvert<1/5$
for all $y\in R({\cal P})^c$.
Moreover, for every $v\in V$,
there exists an integer $m_v$ such that
\[ \lvert m_v-\tilde{\eta}(x)\rvert<\frac{1}{5} \]
for all $x\in X(v,1)$.
Therefore $B_\alpha(\eta)$ is equal to $\sum_{v\in V}-m_v[1_{X(v,1)}]$,
and
\[ \lvert m_v\rvert=\lvert m_v-\tilde{\eta}(\alpha^{h(v)-1}(x))\rvert
\leq\lvert m_v-\tilde{\eta}(x)\rvert+\sum_{k=1}^{h(v)-1}
\lvert\tilde{\eta}(\alpha^{k-1}(x))-\tilde{\eta}(\alpha^k(x))\rvert
<\frac{h(v)}{5}, \]
where $x$ is a point in $X(v,1)$.

Fix $v_0\in V$. Let $a:\{1,2,\dots,h(v_0)\}\rightarrow\R$ be a map
such that $a(h(v_0))=0$, $\lvert m_{v_0}-a(1)\rvert<11/30$ and
$\lvert a(k)-a(k+1)\rvert<11/30$ for all $k=1,2,\dots,h(v_0)-1$.
Define a continuous map $\kappa$ from $X$ to $\R$ by
\[ \kappa(x) \ = \ \begin{cases}
a(k) & x\in X(v_0,k),k=1,2,\dots,h(v_0) \\
\tilde{\eta}(x) & \text{otherwise}. \end{cases} \]
Put $\hat{\kappa}(x)=\kappa(x)+\Z$.
Then $\hat{\kappa}\in C(X,\T)$ belongs to $H(\alpha,\xi)$,
because $7/15-11/30=1/10$.
For $t\in [0,1]$, put $\kappa_t=t\kappa+(1-t)\tilde{\eta}$.
Then it is not hard to see that $\hat{\kappa_t}$ gives a homotopy
from $\eta$ to $\hat{\kappa}$ in $H(\alpha,\xi)$.

At first, let us consider the case that $m_{v_0}$ is positive.
Define $a:\{1,2,\dots,h(v_0)\}\rightarrow\R$ by
$a(1)=a(2)=m_{v_0}+\frac{1}{3}$, $a(3)=m_{v_0}$ and
\[ a(k)=\frac{m_{v_0}(h(v_0)-k)}{h(v_0)-3} \]
for every $k=4,5,\dots,h(v_0)$.
By using this map $a$, define $\kappa\in C(X,\R)$ be as above.
It follows that
$\eta$ and $\hat{\kappa}$ are homotopic in $H(\alpha,\xi)$.
Let $f_1$ be the continuous function on $X\times\T$
defined by
\[ f_1(x,t) \ = \ \begin{cases}
10(t-1/3) &
(x,t)\in X(v_0,2)\times[1/3,1/3+1/10] \\
1 & (x,t)\in X(v_0,2)\times[1/3+1/10,2/3] \\
1-10(t-2/3) & (x,t)\in X(v_0,2)\times[2/3,2/3+1/10] \\
0 & \text{otherwise.} \end{cases} \]
Let $U=X(v_0,2)\times\T$ and
put $f_2=(1_U-f_1)\circ(\alpha\times R_\xi)^{-1}$.
Define a continuous function $g$ on $X\times\T$ by
\[ g(x,t) \ = \ \begin{cases}
-\sqrt{10(t-1/3)(1-10(t-1/3))} &
(x,t)\in X(v_0,2)\times\in [1/3,1/3+1/10] \\
\sqrt{10(t-2/3)(1-10(t-2/3))} &
(x,t)\in X(v_0,2)\times\in [2/3,2/3+1/10] \\
0 & \text{otherwise.} \end{cases} \]
Then it can be verified that $e=gu^*+(f_1+f_2)+ug$ is a projection
and $e$ is equivalent to $1_U$:
in fact, if $h\in C(X\times\T)$ is a function with
\[ h|X(v_0,2)\times[1/3,1/3+1/10]=-1, \ \ \
h|X(v_0,2)\times[2/3,2/3+1/10]=1 \]
and $\lvert h\rvert^2=1$, then the partial isometry
\[ w=h\sqrt{f_1}+\sqrt{1_U-f_1}u^* \]
satisfies $w^*w=e$ and $ww^*=1_U$.
Furthermore $e$ is a subprojection of $e(\alpha,\xi,\hat{\kappa})$
and
\[ e(\alpha,\xi,\hat{\kappa})-e=e(\alpha,\xi,\eta') \]
for some $\eta'\in H(\alpha,\xi)$
with $B_\alpha(\eta')=B_\alpha(\eta)+[1_{X(v_0,1)}]$.
Hence
\begin{align*}
[e(\alpha,\xi,0)]-[e(\alpha,\xi,\eta)]-B_\alpha(\eta)
&= [e(\alpha,\xi,0)]-[e(\alpha,\xi,\hat{\kappa})]-B_\alpha(\eta) \\
&= [e(\alpha,\xi,0)]-[e(\alpha,\xi,\eta')]-[e]-B_\alpha(\eta) \\
&= [e(\alpha,\xi,0)]-[e(\alpha,\xi,\eta')]-B_\alpha(\eta').
\end{align*}
We can repeat the same argument with $\eta'$ in place of $\eta$.
By repeating this $m_{v_0}$ times, we will obtain
\[ [e(\alpha,\xi,0)]-[e(\alpha,\xi,\eta)]-B_\alpha(\eta)
=[e(\alpha,\xi,0)]-[e(\alpha,\xi,\eta')]-B_\alpha(\eta') \]
with $B_\alpha(\eta')=\sum_{v\neq v_0}-m_v[1_{X(v,1)}]$
and $\eta'(x)=0$ for all $x\in X(v_0,1)\cup\dots\cup X(v_0,h(v_0))$.

When $m_{v_0}$ is negative,
there exists $\hat{\kappa}$ homotopic to $\eta$ in $H(\alpha,\xi)$
such that $e(\alpha,\xi,\hat{\kappa})+e=e(\alpha,\xi,\eta')$
for some $\eta'\in H(\alpha,\xi)$
with $B_\alpha(\eta')=B_\alpha(\eta)-[1_{X(v_0,1)}]$.
In a similar fashion to the preceding paragraph,
the same conclusion follows.

By applying the same argument to all towers in $V$,
we have
\[ [e(\alpha,\xi,0)]-[e(\alpha,\xi,\eta)]-B_\alpha(\eta)=0. \]
\end{proof}

\begin{lem}\label{control}
Let $\xa$ be a Cantor minimal system.
Suppose that $\xi_1,\xi_2\in C(X,\R)$ and $f\in C(X,\Z)$
satisfy
\[ \mu(\xi_2)=\mu(\xi_1)+\mu(f) \]
for every $\alpha$-invariant measure $\mu\in M_\alpha$ and
\[ \frac{7}{15}<\xi_i(x)<\frac{8}{15} \]
for all $x\in X$ and $i=1,2$.
Put $\hat{\xi}_i(x)=\xi_i(x)+\Z$ for $i=1,2$.
Then, for any $\varepsilon>0$,
there exists $\eta\in C(X,\T)$ such that
\[ \lvert(\hat{\xi}_1-\hat{\xi}_2)(x)-(\eta-\eta\alpha)(x)\rvert
<\varepsilon \]
for all $x\in X$ and $B_\alpha(\eta)=[f]$.
\end{lem}
\begin{proof}
We may assume $\varepsilon<10^{-1}$.
In the same way as in \cite[Lemma 2.4]{GW},
we can find a Kakutani-Rohlin partition
\[ {\cal P}=\{X(v,k):v\in V,k=1,2,\dots,h(v)\} \]
such that
\[ \frac{1}{h(v)}\left\lvert\sum_{k=1}^{h(v)}
(\xi_1-\xi_2+f)(\alpha^{k-1}(x))\right\rvert<\varepsilon \]
for all $v\in V$ and $x\in X(v,1)$.
Define a real valued continuous function $\kappa$
on $\alpha(R({\cal P}))$ by
\[ \kappa(x)=\sum_{k=1}^{h(v)}
(\xi_1-\xi_2+f)(\alpha^{k-1}(x)) \]
for $x\in X(v,1)$.
Let $\eta$ be a real valued continuous function on $X$ such that
\[ \eta(\alpha^k(x))=\sum_{i=1}^k
(\xi_2-\xi_1)(\alpha^{i-1}(x))+\frac{k}{h(v)}\kappa(x) \]
for all $x\in X(v,1)$ and $k=0,1,\dots,h(v)-1$.
Define $\hat{\eta}(x)=\eta(x)+\Z$.
It is straightforward to check
\[ \lvert(\hat{\xi}_1-\hat{\xi}_2)(x)
-(\hat{\eta}-\hat{\eta}\alpha)(x)\rvert
<\varepsilon \]
for all $x\in X$.
For $x\notin R({\cal P})$, we have
\[ \lvert(\eta-\eta\alpha)(x)\rvert<\frac{1}{15}+\frac{1}{10}, \]
and so the integer nearest to $(\eta-\eta\alpha)(x)$ is zero.
If $x$ belongs to the roof set $R({\cal P})$, then
\[ \lvert(\eta-\eta\alpha)(x)-
\sum_{k=1}^{h(v)}f(\alpha^{1-k}(x))\rvert
<\frac{1}{15}+\frac{1}{10}. \]
Hence we can conclude that $B_\alpha(\eta)$ is equal to $[f]$.
\end{proof}

Now we are ready to prove the main theorem of this section.
Let $\xa$ and $\yb$ be Cantor minimal systems and
let $\xi:X\rightarrow\T$ and $\zeta:Y\rightarrow\T$ be
continuous functions.
Suppose that $\alpha\times R_\xi$ and $\beta\times R_\zeta$
are both minimal.
We denote by $A$ (resp. $B$) the crossed product $C^*$-algebra
arising from $(X\times\T,\alpha\times R_\xi)$
(resp. $(Y\times\T,\beta\times R_\zeta)$).

\begin{thm}\label{KConj1}
The following are equivalent.
\begin{enumerate}
\item $(X\times\T,\alpha\times R_\xi)$ is approximately $K$-conjugate
to $(Y\times\T,\beta\times R_\zeta)$.
\item There exists a unital order isomorphism
$\rho$ from $K_0(B)$ to $K_0(A)$ such that $\rho(K^0\yb)=K^0\xa$.
\end{enumerate}
\end{thm}
\begin{proof}
(1)$\Rightarrow$(2).
This is immediate from the definition of
approximate $K$-conjugacy (Definition \ref{appKconju}).

(2)$\Rightarrow$(1).
By Lemma \ref{perturb},
without loss of generality, we may assume
$\xi(x),\zeta(y)\in(7/15,8/15)$ for all $x\in X$ and $y\in Y$.
Let $e(\alpha,\xi,0)\in A$ and $e(\beta,\zeta,0)\in B$ be
the projections described in the previous section.
We follow the notation used there.
Since $\rho$ is an isomorphism, there are only two possibilities:
\[ \rho([e(\beta,\zeta,0)])-[e(\alpha,\xi,0)]\in K^0\xa \]
or
\[ \rho([e(\beta,\zeta,0)])+[e(\alpha,\xi,0)]\in K^0\xa. \]
Suppose that the latter equality holds.
The dynamical system $(X\times\T,\alpha\times R_\xi)$ is conjugate
to $(X\times\T,\alpha\times R_{-\xi})$ via the mapping
$(x,t)\mapsto (x,-t)$.
This conjugacy induces an isomorphism
between the corresponding $C^*$-algebras, which in turn yields
an isomorphism between the $K_0$-groups.
One can see that $[e(\alpha,\xi,0)]$ is sent to
$[1-e(\alpha,-\xi,0)]$ by this isomorphism.
For this reason,
by replacing $\alpha\times R_\xi$ by $\alpha\times R_{-\xi}$,
we may always assume that there exists $h\in C(X,\Z)$ such that
\[ \rho([e(\beta,\zeta,0)])=[e(\alpha,\xi,0)]+[h]_\alpha \]
in $K_0(A)$.

The restriction of $\rho$ to $K^0\yb$ is
a unital order isomorphism onto $K^0\xa$.
By \cite[Theorem 5.4]{LM} or \cite[Theorem 3.4]{M3},
there exist homeomorphisms $\sigma_n:X\rightarrow Y$
such that
\[ \sigma_n\alpha\sigma_n^{-1}\rightarrow \beta \]
in $\Homeo(Y)$ and
\[ [f\sigma_n]_\alpha=\rho([f]_\beta) \]
for all $f\in C(Y,\Z)$ and $n\in\N$.
We may assume
\[ \lvert\zeta\circ\beta^{-1}\sigma_n(x)
-\zeta\circ\sigma_n\alpha^{-1}(x)\rvert<\frac{1}{n} \]
for all $x\in X$ and $n\in \N$.
As in Remark \ref{states}, we denote the state on the $K_0$-group
arising from an invariant measure $\mu$ by $S_\mu$.
Then, for any $\mu\in M_\alpha$ and $n\in\N$,
we have
\[ \rho^*(S_\mu)([f]_\beta)=S_\mu(\rho([f]_\beta))
=\mu([f\sigma_n]_\alpha)=S_{\sigma_{n*}(\mu)}([f]_\beta) \]
for all $f\in C(Y,\Z)$.
Thus, $\rho^*(S_\mu)=S_{\sigma_{n*}(\mu)}$ on $K^0\yb$,
and so on $K_0(B)$ (see Remark \ref{states}).
Let $\tilde{\xi}\in C(X,\R)$ and $\tilde{\zeta}\in C(Y,\R)$
be the lifts of $\xi$ and $\zeta$ satisfying
\[ \frac{7}{15}<\tilde{\xi}(x)<\frac{8}{15} \ \text{ and } \
\frac{7}{15}<\tilde{\zeta}(y)<\frac{8}{15} \]
for all $x\in X$ and $y\in Y$.
It follows that
\begin{align*}
\mu(\tilde{\zeta}\sigma_n)&= \sigma_{n*}(\mu)(\tilde{\zeta})
=S_{\sigma_{n*}(\mu)}([e(\beta,\zeta,0)]) \\
&= \rho^*(S_\mu)([e(\beta,\zeta,0)])
=S_\mu(\rho([e(\beta,\zeta,0)]) \\
&= S_\mu([e(\alpha,\xi,0)]+[h]_\alpha) \\
&= \mu(\tilde{\xi})+\mu(h)
\end{align*}
for every $\mu\in M_\alpha$ and $n\in\N$.
Now Lemma \ref{control} applies and yields $\eta_n\in C(X,\T)$
such that
\[ \lvert(\xi-\zeta\sigma_n)(x)-(\eta_n-\eta_n\alpha)(x)\rvert
<\frac{1}{n} \]
for all $x\in X$ and $B_\alpha(\eta_n)=[h]$.
Hence it is easy to verify that
\[ (\sigma_n\times R_{\eta_n})(\alpha\times R_\xi)
(\sigma_n\times R_{\eta_n})^{-1}\rightarrow(\beta\times R_\zeta) \]
in $\Homeo(Y\times\T)$.
Furthermore, since
\[ \lvert\zeta\beta^{-1}\sigma_n(x)-
(\xi\alpha^{-1}+\eta_n-\eta_n\alpha^{-1})(x)\rvert<\frac{2}{n} \]
for all $x\in X$, we get the estimate
\[ \lVert f(\beta,\zeta,0)\circ(\sigma_n\times R_{\eta_n})-
f(\alpha,\xi,-\eta_n)\rVert<\frac{20}{n}. \]
See the discussion before Proposition \ref{RieffelI}
for the definition of $f(\cdot,\cdot,\cdot)$.
It is easy to see
$g_0\circ(\sigma_n\times R_{\eta_n})=g_{-\eta_n}$.
Therefore, the asymptotic morphism $\{\psi_n\}:B\rightarrow A$
associated with $\sigma_n\times R_{\eta_n}$ satisfies
\[ \lim_{n\rightarrow\infty}
\lVert\psi_n(e(\beta,\zeta,0))-e(\alpha,\xi,-\eta_n)\rVert=0, \]
and Lemma \ref{Bott} yields
\begin{align*}
[e(\alpha,\xi,-\eta_n)]&= [e(\alpha,\xi,0)]+B_\alpha(\eta_n) \\
&= [e(\alpha,\xi,0)]+[h]_\alpha=\rho([e(\beta,\zeta,0)]).
\end{align*}
For every clopen set $U\subset Y$, we know
\[ \lim_{n\rightarrow\infty}
\lVert\psi_n(1_U)-1_U\circ\sigma_n\rVert=0 \]
and $[1_U\circ\sigma_n]_\alpha=\rho([1_U]_\beta)$.
It follows that $\{\psi_n\}$ induces
$\rho:K_0(B)\rightarrow K_0(A)$.

For every clopen set $U\subset Y$, we know
\[ \lim_{n\rightarrow\infty}
\lVert\psi_n(z1_U)-z1_U\circ(\sigma_n\times R_{\eta_n})\rVert=0, \]
where $z$ is a unitary defined by $z(y,t)=e^{2\pi\sqrt{-1}t}$.
It is clear that
\[ [1_{U^c}+z1_U\circ(\sigma_n\times R_{\eta_n})]
=[1_{U^c}+z(1_U\circ\sigma_n)] \]
in $K_1(A)$.
Since $\{\psi_n\}$ approximately carries the implementing unitary
of $\beta\times R_\zeta$ to that of $\alpha\times R_\xi$,
we can conclude that $\{\psi_n\}$ induces an isomorphism
between $K_1(B)$ and $K_1(A)$.

Consequently $\{\psi_n\}$ induces an isomorphism between $K$-groups.
Similarly, we can construct an asymptotic morphism
from $A$ to $B$ which induces $\rho^{-1}$ between their $K$-groups.
\end{proof}

\begin{thm}\label{KConj2}
Suppose that $\alpha\times R_\xi$ and $\beta\times R_\zeta$
are minimal and rigid. Then the following are equivalent.
\begin{enumerate}
\item $(X\times\T,\alpha\times R_\xi)$ is approximately $K$-conjugate
to $(Y\times\T,\beta\times R_\zeta).$
\item There exists a unital order isomorphism
$\rho$ from $K_0(B)$ to $K_0(A)$ such that $\rho(K^0\yb)=K^0\xa.$
\item $(X\times\T,\alpha\times R_\xi)$ is approximately flip $K$-conjugate
to $(Y\times\T,\beta\times R_\zeta).$
\item $\alpha\times R_\xi$ and $\beta\times R_\zeta$ are
$C^*$-strongly approximately flip conjugate.
\item There is $\theta\in KL(A,B)$ which gives an order and
unit preserving isomorphism from $(K_0(A),K_0(A)_+,[1_A],K_1(A))$
onto $(K_0(B),K_0(B)_+,[1_B],K_1(B))$ and an isomorphism
$\chi:C(X\times\T)\to C(Y\times\T)$ such that
$$
[j_{\alpha\times R_{\xi}}]\times\theta=
[j_{\beta\times R_\zeta}\circ \chi]
$$
in $KL(C(X\times\T), B).$
\end{enumerate}
\end{thm}
\begin{proof}
We have seen (1)$\Leftrightarrow$(2) in Theorem \ref{KConj1}.

(1)$\Rightarrow$(3) is obvious.

When both $\alpha\times R_\xi$ and $\beta\times R_\zeta$ are
rigid, $A$ and $B$ have tracial rank zero by Theorem \ref{Tad1}.
Thus,(3)$\Rightarrow$(4) follows from \cite[Theorem 5.4]{Lnd}.

(4)$\Rightarrow$(5) follows immediately from \cite[Theorem 3.9]{Lnd}.

(5)$\Rightarrow$(2).
Suppose that $\theta$ induces an order and unit preserving isomorphism
$\Gamma(\theta)$ from $(K_0(A),K_0(A)_+,[1_A], K_1(A))$ to
$(K_0(B), K_0(B)_+, [1_B], K_1(B)).$ Suppose that
there is $\chi: C(X\times \T)\to C(Y\times \T)$
such that
$[j_{\alpha\times R_\xi}]\times\theta=
[j_{\beta\times R_\zeta}\circ \chi].$
This implies that $\Gamma(\theta)$ gives an isomorphism
from $K^0(X,\alpha)$ onto $K^0(Y, \beta).$ So (2) holds.
\end{proof}

\section{Non-orientation preserving case}

Let $\xa$ be a Cantor minimal system and
let $\phi:X\rightarrow\Homeo(\T)$ be a continuous map.
In this section, we would like to consider
the case that $\alpha\times\phi$ is not orientation preserving,
that is, $[o(\phi)]$ is not zero in $K^0\xa/2K^0\xa$.
As was seen in the discussion before Lemma \ref{Kofnop},
the skew product extension $(X\times\Z_2,\alpha\times o(\phi))$ is
a Cantor minimal system.
This system will play an important role
when we study $\alpha\times\phi$.

Define a continuous map
$\widetilde{\phi}:X\times\Z_2\rightarrow\Homeo(\T)^+$ by
\[ \widetilde{\phi}_{(x,k)}=\lambda^{k+o(\phi)(x)}\phi_x\lambda^k, \]
where $\lambda$ is given by $\lambda(t)=-t$ for $t\in\T$.
Let $\pi$ be the projection from $X\times\Z_2$ to
the first coordinate.

\begin{lem}\label{untwist}
As a $\Homeo(\T)$-valued cocycle on the Cantor minimal system
$(X\times\Z_2,\alpha\times o(\phi))$,
$\phi\pi$ is cohomologous to $\widetilde{\phi}$.
In particular, $\phi\pi$ is orientation preserving with respect to
the minimal homeomorphism $\alpha\times o(\phi)$.
\end{lem}
\begin{proof}
Put $\omega_{(x.k)}=\lambda^k$. Then
\[ \widetilde{\phi}_{(x,k)}\circ\omega_{(x,k)}
=\lambda^{k+o(\phi)(x)}\circ\phi_x
=\omega_{(\alpha\times o(\phi))(x,k)}\circ\phi_{\pi(x,k)} \]
implies that they are cohomologous.
\end{proof}

We remark that the following diagram of factor maps is commutative:
\[ \begin{CD}
(X\times\T,\alpha\times\phi) @<\pi\times\id<<
(X\times\Z_2\times\T,\alpha\times o(\phi)\times\phi\pi) \\
@VVV @VVV \\
\xa @<<\pi< (X\times\Z_2,\alpha\times o(\phi)).
\end{CD} \]

From Lemma \ref{prepaA} and \ref{prepaB}, there exists $\xi\in C(X,\T)$
such that $r(R_{\xi(x)}\lambda^{o(\phi)(x)}\phi_x)\neq0$ for all $x\in X$.
By applying Lemma \ref{perturbIII} to $-\xi$ and $o(\phi)$,
we obtain $\eta\in C(X,\T)$ such that
\[ \lvert\eta(x)-(-1)^{o(\phi)(x)}\eta\alpha(x)-\xi(x)\rvert
<\varepsilon, \]
where $\varepsilon$ is sufficiently small so that it implies
\begin{align*}
0&\neq r(R_{\eta(x)-(-1)^{o(\phi)(x)}\eta\alpha(x)}
\lambda^{o(\phi)(x)}\phi_x) \\
&= r(R_{(-1)^{o(\phi)(x)}\eta\alpha(x)}^{-1}
\lambda^{o(\phi)(x)}\phi_xR_{\eta(x)}) \\
&= r(\lambda^{o(\phi)(x)}R_{\eta\alpha(x)}^{-1}\phi_xR_{\eta(x)})
\end{align*}
for all $x\in X$.
Therefore, by perturbing $\phi$ by $R_\eta$,
we may assume
\[ r(\widetilde{\phi}_{(x,k)})
=(-1)^kr(\lambda^{o(\phi)(x)}\phi_x)\neq0 \]
for all $(x,k)\in X\times\Z_2$.

Let $A$ denote the crossed product $C^*$-algebra arising
from $(X\times\T,\alpha\times\phi)$.
We write the implementing unitary by $u$.
Define an automorphism $\theta\in\Aut(A)$ of order two by
\[ \theta(f)=f \]
for all $f\in C(X\times\T)$ and
\[ \theta(u)=ug, \]
where $g\in C(X\times\T)$ is given by $g(x,t)=(-1)^{o(\phi)(x)}$.
By Lemma \ref{perturb}, it can be seen that
$\theta$ is approximately inner, and so
it induces the identity on the $K$-group.

\begin{prop}
In the situation above, the crossed product $C^*$-algebra
arising from the dynamical system
$(X\times\Z_2\times\T,\alpha\times o(\phi)\times \phi\pi)$
is isomorphic to $A\rtimes_\theta\Z_2$.
\end{prop}
\begin{proof}
We write the implementing unitary in $A$ by $u_A$.
Let us denote
$C^*(X\times\Z_2\times\T,\alpha\times o(\phi)\times \phi\pi)$
by $B$ and the implementing unitary in $B$ by $u_B$.
Let $v$ denote the unitary which implements $\theta$.

We would like to define a homomorphism $\Phi$
from $A\rtimes\Z_2$ to $B$.
The $C^*$-algebra $A$ can be naturally embedded into $B$
via the factor map $\pi\times\id$
from $X\times\Z_2\times\T$ to $X\times\T$.
Let $\Phi|A$ be this embedding.
Define a continuous function $h\in C(X\times\Z_2\times\T)$ by
\[ h(x,k,t)=(-1)^k. \]
Then,
\[ h\circ(\alpha\times o(\phi)\times\phi\pi)=
h(g\circ(\pi\times\id)). \]
Put $\Phi(v)=h$.
For $f\in C(X\times\T)$, we have
\[ h\Phi(f)h=h(f\circ(\pi\times\id))h=f\circ(\pi\times\id)
=\Phi(\theta(f)). \]
Besides,
\[ h\Phi(u_A)h=hu_Bh=u_B(g\circ(\pi\times\id))
=\Phi(\theta(u_A)). \]
It follows that $\Phi$ is a well-defined homomorphism.
It is not hard to see that $\Phi$ is an isomorphism.
\end{proof}

We freely use the identification of the two $C^*$-algebras established
in the proposition above.
By Lemma \ref{untwist} and \ref{Kofop}, we know that
\[ K_0(A\rtimes\Z_2)\cong\Z\oplus K^0(X\times\Z_2,\alpha\times o(\phi)) \]
and
\[ K_1(A\rtimes\Z_2)\cong\Z\oplus K^0(X\times\Z_2,\alpha\times o(\phi)). \]
On the crossed product
$C^*(X\times\Z_2\times\T,\alpha\times o(\phi)\times\phi\pi)$,
the dual action $\hat{\theta}$ is given by
\[ \hat{\theta}(f)(x,k,t)=f(x,k+1,t) \]
for $f\in C(X\times\Z_2\times\T)$ and
\[ \hat{\theta}(u)=u, \]
where $u$ is the implementing unitary.
Define a homeomorphism $\gamma\in\Homeo(X\times\Z_2)$ by
\[ \gamma(x,k)=(x,k+1). \]
Then $\gamma\times\id$ commutes with $\alpha\times o(\phi)\times\phi\pi$
and $\hat{\theta}$ on $C(X\times\Z_2\times\T)$ is
induced by $\gamma\times\id$.
Let us consider the induced action $\hat{\theta}_*$ on the $K$-groups.
Evidently,
$\hat{\theta}_{*0}$ on $K^0(X\times\Z_2,\alpha\times o(\phi))$ is given by
\[ [f]\mapsto[f\gamma], \]
and
$\hat{\theta}_{*1}$ on $K^0(X\times\Z_2,\alpha\times o(\phi))$ is given by
\[ [f]\mapsto[-f\gamma]. \]
Of course, $\hat{\theta}_*([u])=[u]$.
Hence it remains to know the image of the generalized Rieffel projection.
Let
\[ e=e(\alpha\times o(\phi),\widetilde{\phi},0)
=g_0u^*+f(\alpha\times o(\phi),\widetilde{\phi},0)+ug_0 \]
be the projection of
$C^*(X\times\Z_2\times\T,\alpha\times o(\phi)\times\widetilde{\phi})$
as in Proposition \ref{RieffelII}.
This is well-defined, because $r(\widetilde{\phi}_{(x,k)})$ is not zero.
By the map
\[ (x,k,t)\mapsto(x,k+1,\lambda(t)), \]
$f(\alpha\times o(\phi),\widetilde{\phi},0)$ is carried to
a function supported on
\[ \left\{(x,k,t):
t\in[-\widetilde{\phi}_{(\alpha\times o(\phi))^{-1}(x,k+1)}(c),0]\right\} \]
and $g_0$ is carried to a function supported on $X\times\Z_2\times[-c,0]$.
Note that
\begin{align*}
& -\widetilde{\phi}_{(\alpha\times o(\phi))^{-1}(x,k+1)}(c) \\
&= -\lambda^{k+1}\phi_{\alpha^{-1}(x)}
\lambda^{k+1+o(\phi)(\alpha^{-1}(x))}(c) \\
&= -\lambda\widetilde{\phi}_{(\alpha\times o(\phi))^{-1}(x,k)}\lambda(c) \\
&= \widetilde{\phi}_{(\alpha\times o(\phi))^{-1}(x,k)}(-c).
\end{align*}
Hence, under the identification of
$C^*(X\times\Z_2\times\T,\alpha\times o(\phi)\times\widetilde{\phi})$
with
$C^*(X\times\Z_2\times\T,\alpha\times o(\phi)\times\phi\pi)$,
we have
\begin{align*}
& \hat{\theta}(e)=\hat{\theta}(e(\alpha\times o(\phi),\widetilde{\phi},0)) \\
&= g_{-c}u^*+1-f(\alpha\times o(\phi),\widetilde{\phi},-c)+ug_{-c} \\
&= 1-(-g_{-c}u^*+f(\alpha\times o(\phi),\widetilde{\phi},-c)-ug_{-c}).
\end{align*}
By the remark following Proposition \ref{RieffelII},
this is homotopic to
\[ 1-e(\alpha\times o(\phi),\widetilde{\phi},-c), \]
and moreover it is homotopic to
\[ 1-e(\alpha\times o(\phi),\widetilde{\phi},0) \]
from Proposition \ref{RieffelII}. Thus, we get
\[ [\hat{\theta}(e)]=1-[e] \]
in $K_0(A\rtimes\Z_2)$.
In particular, $\hat{\theta}$ is not approximately inner.

Next, we would like to consider the map between $K$-groups induced
from the inclusion $\iota:A\hookrightarrow A\rtimes\Z_2$.
On the $K_0$-group, that is clearly given by
\[ \iota_{*0}([f])=
(0,[f\circ\pi])\in\Z\oplus K^0(X\times\Z_2,\alpha\times o(\phi))\cong
K_0(A\rtimes\Z_2) \]
for $[f]\in K^0\xa\cong K_0(A)$.
On the $K_1$-group, for $[f]\in \Coker(\id-\alpha_\phi^*)$,
we can see that
\[ \iota_{*1}([f])=
(0,[\delta(f)])\in\Z\oplus K^0(X\times\Z_2,\alpha\times o(\phi))\cong
K_1(A\rtimes\Z_2), \]
where $\delta(f)(x,k)=(-1)^kf(x)$ for $(x,k)\in X\times\Z_2$.
Under the identification
\[ \Coker(\id-\alpha_\phi^*)\cong
K^0(X\times\Z_2,\alpha\times o(\phi))/K^0\xa \]
established in the discussion before Lemma \ref{Kofnop},
this map is also described by
\[ \iota_{*1}([f]+K^0\xa)=(0,[f-f\gamma]) \]
for $f\in C(X\times\Z_2,\Z)$.
\bigskip

We now consider minimality and rigidity of $\alpha\times\phi$.

\begin{lem}\label{nopminimal}
Let $\xa$ be a Cantor minimal system and
let $\phi:X\rightarrow\Homeo(\T)$ be a continuous map.
Suppose that $\alpha\times\phi$ is not orientation preserving.
Then $\alpha\times\phi$ is minimal if and only if
$\alpha\times o(\phi)\times\widetilde{\phi}$ is minimal.
\end{lem}
\begin{proof}
We follow the notation used in the discussion above.
By Lemma \ref{untwist}, we may replace
$\alpha\times o(\phi)\times\widetilde{\phi}$ by
$\alpha\times o(\phi)\times\phi\pi$.

Suppose that $\alpha\times o(\phi)\times\phi\pi$ is minimal.
If $E\subset X\times\T$ is a closed $\alpha\times\phi$-invariant
subset, then $(\pi\times\id)^{-1}(E)$ is a closed
$\alpha\times o(\phi)\times\phi\pi$-invariant subset of
$X\times\Z_2\times\T$.
It follows that $(\pi\times\id)^{-1}(E)$ is empty or
the whole $X\times\Z_2\times\T$.
Namely $E$ is empty or the whole $X\times\T$.

Let us prove the converse.
Assume that $\alpha\times\phi$ is minimal.
Let $E\subset X\times\Z_2\times\T$ be a minimal subset
of $\alpha\times o(\phi)\times\phi\pi$.
Since $(\pi\times\id)(E)$ is a closed $\alpha\times\phi$-invariant
subset, it must be equal to $X\times\T$.
If $E=\gamma(E)$, then $E=X\times\Z_2\times\T$ and
we have nothing to do.
Suppose that $E\cap\gamma(E)$ is empty and
$E\cup\gamma(E)=X\times\Z_2\times\T$.
Thus $E$ is a clopen subset.
Hence there exists a continuous function
$\chi:X\rightarrow\Z_2$ such that
\[ E=\{(x,\chi(x),t):x\in X,t\in\T\}. \]
It follows that $\chi+o(\phi)=\chi\alpha$,
which contradicts $[o(\phi)]\neq0$ in $K^0\xa/2K^0\xa$.
\end{proof}

\begin{lem}
Let $\xa$ be a Cantor minimal system and
let $\phi:X\rightarrow\Homeo(\T)$ be a continuous map.
Suppose that $\alpha\times\phi$ is not orientation preserving.
\begin{enumerate}
\item If $\alpha\times o(\phi)\times\widetilde{\phi}$ is rigid,
then $\alpha\times\phi$ is rigid.
\item If $\phi$ takes its values in $\Isom(\T)$ and
$\alpha\times\phi$ is rigid, then
$\alpha\times o(\phi)\times\widetilde{\phi}$ is rigid.
\end{enumerate}
\end{lem}
\begin{proof}
We follow the notation used in the discussion above.
By Lemma \ref{untwist}, we may replace
$\alpha\times o(\phi)\times\widetilde{\phi}$ by
$\alpha\times o(\phi)\times\phi\pi$.
Let $F$ (resp. $\widetilde{F}$) denote the canonical factor map
from $(X\times\T,\alpha\times\phi)$ to $\xa$ (resp. from
$(X\times\Z_2\times\T,\alpha\times o(\phi)\times\phi\pi)$
to $(X\times\Z_2,\alpha\times o(\phi))$.

(1) Suppose that there exist two distinct ergodic measures
$\nu_1$ and $\nu_2$ for $(X\times\T,\alpha\times\phi)$
such that $F_*(\nu_1)=F_*(\nu_2)$.
Let $\tilde{\nu}_i$ be
an $\alpha\times o(\phi)\times\phi\pi$-invariant measure
such that $(\pi\times\id)_*(\tilde{\nu}_i)=\nu_i$.
Of course, $\tilde{\nu}_1\neq\tilde{\nu}_2$,
because of
\[ (\pi\times\id)_*(\tilde{\nu}_1)=\nu_1\neq
\nu_2=(\pi\times\id)_*(\tilde{\nu}_2). \]
By replacing $\tilde{\nu}_i$ by
\[ \frac{1}{2}(\tilde{\nu}_i+(\gamma\times\id)_*(\tilde{\nu}_i)), \]
we may assume that $\tilde{\nu}_i$ is $\gamma\times\id$-invariant.
It follows that $\widetilde{F}_*(\tilde{\nu}_i)$ is invariant
under $\gamma$.
Together with
\[ \pi_*\widetilde{F}_*(\tilde{\nu}_1)=
F_*(\pi\times\id)_*(\tilde{\nu}_1)=
F_*(\nu_1)=F_*(\nu_2)=
F_*(\pi\times\id)_*(\tilde{\nu}_2)=
\pi_*\widetilde{F}_*(\tilde{\nu}_2), \]
we have $\widetilde{F}_*(\tilde{\nu}_1)=
\widetilde{F}_*(\tilde{\nu}_2)$.
Therefore $\alpha\times o(\phi)\times\phi\pi$ is not rigid.

(2) Assume that $\alpha\times o(\phi)\times\phi\pi$ is not rigid.
There exists an ergodic measure $\mu$ for
$(X\times\Z_2,\alpha\times o(\phi))$ such that
$\widetilde{F}_*^{-1}(\mu)$ is not a singleton.
By assumption, $\widetilde{\phi}$ takes its values
in the rotation group.
It follows from Lemma \ref{whenrigid} and its proof that
there exist uncountably many ergodic measures for
$(X\times\Z_2\times\T,\alpha\times o(\phi)\times\phi\pi)$
in $\widetilde{F}_*^{-1}(\mu)$.
In particular, we can find two distinct ergodic measures
$\nu_1,\nu_2\in \widetilde{F}_*^{-1}(\mu)$ such that
$(\gamma\times\id)_*(\nu_1)\neq\nu_2$.
Hence, it is easily verified that
\[ (\pi\times\id)_*(\nu_1)\neq(\pi\times\id)_*(\nu_2) \]
in $M_{\alpha\times\phi}$.
But, we have
\[ F_*(\pi\times\id)_*(\nu_1)=
\pi_*\widetilde{F}_*(\nu_1)=\pi_*(\mu)=
\pi_*\widetilde{F}_*(\nu_2)=
F_*(\pi\times\id)_*(\nu_2), \]
and so $\alpha\times\phi$ is not rigid.
\end{proof}
\bigskip

Now we consider cocycles with values in $\Isom(\T)$.
Let $\xa$ be a Cantor minimal system and
let $\phi:X\rightarrow\Isom(\T)$ be a continuous map.
There exists $\xi\in C(X,\T)$ such that
$\phi_x=\lambda^{o(\phi)(x)}R_{\xi(x)}$ for all $x\in X$.
Suppose that $\alpha\times\phi$ is not orientation preserving and
not minimal.
By Lemma \ref{nopminimal},
$\alpha\times o(\phi)\times\widetilde{\phi}$ is not minimal,
where $\widetilde{\phi}$ is given by
\[ \widetilde{\phi}_{(x,k)}=
\lambda^{k+o(\phi)(x)}\phi_x\lambda^k=
\lambda^kR_{\xi(x)}\lambda^k=
\begin{cases}R_{\xi(x)} & k=0 \\
R_{-\xi(x)} & k=1. \end{cases} \]
It is convenient to introduce
$\delta:C(X,\T)\rightarrow C(X\times\Z_2,\T)$ defined by
$\delta(\eta)(x,k)=(-1)^k\eta(x)$.
It follows from Lemma \ref{opminimal} that
there exist $n\in\N$ and $\zeta\in C(X\times\Z_2,\Z)$ such that
\[ n\delta(\xi)=\zeta-\zeta\circ(\alpha\times o(\phi))^{-1}. \]
We also have
\[ -n\delta(\xi)=n\delta(\xi)\circ\gamma=
(\zeta-\zeta\circ(\alpha\times o(\phi))^{-1})\circ\gamma=
\zeta\gamma-\zeta\gamma\circ(\alpha\times o(\phi))^{-1}. \]
Since $\alpha\times o(\phi)$ is minimal on $X\times\Z_2$,
$\zeta+\zeta\gamma$ must equal a constant function.
We can adjust $\zeta$ by a constant function
so that $\zeta+\zeta\gamma$ is equal to zero.
Thus there exists $\eta\in C(X,\T)$ such that
\[ n\delta(\xi)=
\delta(\eta)-\delta(\eta)\circ(\alpha\times o(\phi))^{-1}. \]
Combining this with
$\delta\circ\alpha_\phi^*=(\alpha\times o(\phi))^*\circ\delta$,
we obtain
\[ n\xi=\eta-\alpha_\phi^*(\eta). \]

\begin{lem}\label{Isomminimal}
Let $\xa$ be a Cantor minimal system and
let $\phi_x=\lambda^{o(\phi)(x)}R_{\xi(x)}$ be a cocycle
with values in $\Isom(\T)$.
If $\alpha\times\phi$ is not orientation preserving and
not minimal, then there exist $n\in\N$ and $\eta\in C(X,\T)$
such that
\[ n\xi=\eta-\alpha_\phi^*(\eta). \]
Moreover, every minimal subset of $\alpha\times\phi$ is given by
\[ E_s=\{(x,t):
nt=\alpha_\phi^*(\eta)(x)+s \ \text{ or } \
nt=\alpha_\phi^*(\eta)(x)-s\} \]
for some $s\in\T$.
\end{lem}
\begin{proof}
The first part follows the discussion above.
Let us consider the latter part.
By Lemma \ref{opminimal}, every minimal set of
$\alpha\times o(\phi)\times\widetilde{\phi}$ is given by
\[ \{(x,k,t):nt=
(-1)^{k+o(\phi)(\alpha^{-1}(x))}\eta\alpha^{-1}(x)+s\}. \]
This closed set is carried to
\[ E_s=\{(x,t):
nt=\alpha_\phi^*(\eta)(x)+s \ \text{ or } \
nt=\alpha_\phi^*(\eta)(x)-s\} \]
by the factor map from $X\times\Z_2\times\T$ to $X\times\T$.
Therefore $E_s$ is a minimal set of $\alpha\times\phi$.
\end{proof}

\section{Examples}

\begin{exm}\label{Furstenberg}
Let $\theta\in\T$ be an irrational number and
let $\xi:\T\rightarrow\T$ be a continuous map.
A homeomorphism
\[ \gamma:(s,t)\mapsto(s+\theta,t+\xi(s)) \]
on $\T^2$ is called a Furstenberg transformation.
In \cite{OP}, the crossed product $C^*$-algebra arising from
$(\T^2,\gamma)$ is studied.
We would like to replace the irrational rotation
with a Cantor minimal system and
construct an almost one to one extension of $(\T^2,\gamma)$
as follows.
Let $\phi$ be a Denjoy homeomorphism on $\T$
with $r(\phi)=\theta$ as in Remark \ref{rigidbutnotminimal}.
The homeomorphism $\phi$ has
the unique invariant nontrivial closed subset $X$.
Let $\alpha$ be the restriction of $\phi$ to $X$.
Then $\xa$ is a Cantor minimal system and
there exists an almost one to one factor map $\pi$ from $\xa$
to the irrational rotation $(\T,R_\theta)$
(see \cite{PSS} for details).
Both $\xa$ and $(\T,R_\theta)$ are uniquely ergodic.
It is easy to see that $\pi\times\id:X\times\T\rightarrow\T^2$
satisfies
\[ (\pi\times\id)\circ(\alpha\times R_{\xi\pi})=
\gamma\circ(\pi\times\id), \]
that is, $\pi\times\id$ is a factor map.
One can check that if $\gamma$ is minimal, then
$\alpha\times R_{\xi\pi}$ is also minimal.
The factor map $\pi$ induces a Borel isomorphism
between $X$ and $S^1$.
Hence $\alpha\times R_{\xi\pi}$ is rigid if and only if
$\gamma$ is uniquely ergodic.
There are some known criterions for unique ergodicity of $\gamma$.
For example, it was proved in \cite[Theorem 2.1]{F} that
if $\xi$ is a Lipschitz function and its degree is not zero
then $\gamma$ is uniquely ergodic.
On the other hand
it is known that there exist $\theta\in\T$ and
$\xi:\T\rightarrow\T$ such that $\gamma$ is minimal but not
uniquely ergodic (see \cite[p.585]{F} for instance).
In this case
$\alpha\times R_{\xi\pi}$ is minimal but not rigid.
\end{exm}

Now let us consider the example of Putnam which was presented by
N. C. Phillips in \cite{Ph3}.

\begin{exm}\label{Exadd1}
For any $\theta\in \R\setminus \Q$ let $g_{\theta}$ be a minimal
homeomorphism of a Cantor set $X_{\theta}\subset \T$ obtained from
a Denjoy homeomorphism $g_{\theta}^0: \T\to \T$ as in \cite{Ph3}.
Choose $g_{\theta}^0$ to have rotation number $\theta$ and
such that the unique minimal set $X_{\theta}\subset \T$
has the property that
the image of $\T\setminus X_{\theta}$ under the semiconjugation to
$R_{\theta}$ is a single orbit of $R_{\theta}.$

Now let $\theta_1, \theta_2\in \R\setminus \Q$ be irrational
numbers such that $1, \theta_1, \theta_2$ are $\Q$-linearly
independent. Consider two systems $(X_{\theta_1}\times \T,
g_{\theta_1}\times R_{\theta_2})$ and
$(X_{\theta_2}\times \T,g_{\theta_2}\times R_{\theta_1}).$ Then
both are minimal.  Let $A=C^*(X_{\theta_1}\times \T,
g_{\theta_1}\times R_{\theta_2})$ and $B=C^*(X_{\theta_2}\times
\T, g_{\theta_2}\times R_{\theta_1}).$
It follows from \cite[Proposition 1.12]{Ph3} that
$$
K_1(A)\cong K_1(B)=\Z^3
$$
and
$$
(K_0(A),K_0(A)_+, [1_A])\cong (\Z+\theta_1\Z+\theta_2\Z,
(\Z+\theta_1\Z+\theta_2\Z)_+,1)\cong (K_0(B), K_0(B)_+,[1_B]),
$$
where $\Z+\theta_1\Z+\theta_2\Z\subset \R.$
In particular, both systems are rigid. Moreover $A$ and $B$ have
tracial rank zero and they are isomorphic by the classification
theorem in \cite{L3}.  However, there is no order isomorphism
between $\Z+\theta_1\Z$ and $\Z+\theta_2\Z.$ It follows from
Theorem \ref{KConj2} that they are not approximately $K$-conjugate.
On the other hand, by \cite[Corollary 4.10]{M3},
they are weakly approximately conjugate.

When we choose another continuous function $\xi:X_{\theta_1}\to\T$,
two systems $g_{\theta_1}\times R_{\theta_2}$ and
$g_{\theta_1}\times R_{\xi}$ may not be conjugate.
But, if the integral value of $\xi$ is equal to $\theta_2$, then
we can conclude that they are approximately $K$-conjugate
by Theorem \ref{KConj1}.

\end{exm}

%
%
%

\begin{exm}
We would like to construct a non-orientation preserving
minimal homeomorphism on $X\times\T$ concretely.
Let $\xa$ be an odometer system of type $3^\infty$.
It is well-known that $K^0\xa$ is isomorphic to $\Z[1/3]$.
We regard $X$ as a projective limit of $\Z_{3^n}$ and
denote the canonical projection from $X$ to $\Z_{3^n}$
by $\pi_n$.
Notice that $\pi_n(\alpha(x))=\pi_n(x)+1$ for all $x\in X$,
where the addition is understood modulo $3^n$.
Let $x_0$ be the point of $X$ such that
$\pi_n(x_0)=0$ for all $n\in\N$.
We will construct a continuous map $\phi:X\rightarrow\Isom(\T)$
of the form $\phi_x=\lambda R_{\xi(x)}$ so that
$\alpha\times\phi$ is a minimal homeomorphism on $X\times\T$.
Note that $o(\phi)(x)=1$ and $[o(\phi)]\neq0$
in $K^0\xa/2K^0\xa\cong\Z_2$.
By Lemma \ref{Isomminimal}, if the closure of
\[ \{(\alpha\times\phi)^m(x_0,0):m\in\N\} \]
contains $\{x_0\}\times\T$, then we can deduce the minimality of
$\alpha\times\phi$.
Let $\{t_n\}_{n\in\N}$ be a dense sequence of $\T$.
Since $\alpha^{3^n}(x_0)\rightarrow x_0$ as $n\rightarrow\infty$,
it suffices to construct $\phi$ so that
$(\alpha\times\phi)^{3^n}(x_0,0)=(\alpha^{3^n}(x_0),-t_n)$.

Let $s_n\in (-2^{-1},2^{-1}]$ be the real number satisfying
$s_n+\Z=t_n-t_{n-1}$, where we put $t_0=0$.
We define a map $\xi_n:X\rightarrow\T$ by
\[ \xi_n(x)=\begin{cases}
0 & \pi_n(x)=0,1,\dots,3^{n-1}-1 \\
\frac{(-1)^k}{2\cdot3^{n-1}}s_n+\Z & \text{otherwise}. \end{cases}\]
For all $n,m\in\N$, it is not hard to see that
\[ \sum_{k=0}^{3^n-1}(-1)^k\xi_m(\alpha^k(x_0))=
\begin{cases} 0 & n<m \\ s_m+\Z & n\geq m. \end{cases} \]
Since $\lvert\xi_n(x)\rvert<3^{-n}$ for every $x\in X$,
\[ \xi(x)=\sum_{n=1}^\infty\xi_n(x) \]
exists and is continuous on $X$.
Put $\phi_x=\lambda R_{\xi(x)}$ for all $x\in X$.
Then
\[ \sum_{k=0}^{3^n-1}(-1)^k\xi(\alpha^k(x_0))=
\sum_{i=1}^ns_i+\Z=t_n \]
implies
\[ (\alpha\times\phi)^{3^n}(x_0,0)=(\alpha^{3^n}(x_0),\lambda(t_n))
=(\alpha^{3^n}(x_0),-t_n) \]
for all $n\in\N$.
It follows that $\alpha\times\phi$ is minimal.
\end{exm}

\begin{exm}
We would like to construct a cocycle with values in $\Homeo^+(\T)$
which is not cohomologous to a cocycle with values
in the rotation group.
It is useful to introduce a complete metric $d(\cdot,\cdot)$
of $\Homeo^+(\T)$ defined by
\[ d(\phi,\psi)=\max_{t\in\T}
\left\{\lvert\phi(t)-\psi(t)\rvert,
\lvert\phi^{-1}(t)-\psi^{-1}(t)\rvert\right\}. \]
In the argument below, we use the following facts.
\begin{description}
\item[Fact (a).] For any $s,t\in\T$ and $\varepsilon>0$,
there exists $\rho\in\Homeo^+(\T)$ such that
$\lvert\rho(s)-s\rvert<\varepsilon$,
$\lvert\rho(t)-s\rvert<\varepsilon$ and
$\rho$ is conjugate to an irrational rotation.
\item[Fact (b).] $\Homeo^+(\T)$ is arcwise connected.
\end{description}
Let us construct two sequences of natural numbers $m_n$ and $l_n$,
and a sequence of maps $\phi_n:\Z_{m_n}\rightarrow\Homeo^+(\T)$
inductively so that the following conditions are satisfied.
It is convenient to view $\phi_n$ as a periodic map from $\Z$.
\begin{enumerate}
\item $m_{n-1}$ divides $m_n$ and
$l_n$ is not greater than $m_n/m_{n-1}$.
\item $\psi_n=\phi_n(m_n-1)\dots\phi_n(1)\phi_n(0)$ is conjugate to
an irrational rotation.
\item Both $\lvert\psi_n(0)\rvert$ and $\lvert\psi_n(1/2)\rvert$ are
less than $1/n$.
\item For every $t\in\T$, $\{\psi_{n-1}^k(t):k=1,2,\dots,l_n\}$
is $1/n$-dense in $\T$.
\item For every $k=0,1,\dots,l_nm_{n-1}-1$,
$\phi_n(k)=\phi_{n-1}(k)$.
\item For every $k\in\Z_{m_n}$,
$d(\phi_n(k),\phi_{n-1}(k))$ is less than $2^{-n}$.
\end{enumerate}
If these conditions are achieved,
then we can finish the proof as follows.
Let $\xa$ be the odometer system of type $\{m_n\}_n$.
Namely, $X$ is the projective limit of $\Z_{m_n}$ and
there exists a natural projection $\pi_n:X\rightarrow\Z_{m_n}$
such that $\pi_n(\alpha(x))=\pi_n(x)+1$,
where the addition is understood modulo $m_n$.
Let $x_0$ be the point of $X$ such that $\pi_n(x_0)=0$ for all $n\in\N$.
From (6),
\[ \phi_x=\lim_{n\rightarrow\infty}\phi_n(\pi_n(x))\in\Homeo^+(\T) \]
exists for all $x\in X$ and
$\phi$ is a continuous map from $X$ to $\Homeo^+(\T)$.
By (4) and (5), for all $n\in\N$ and $t\in\T$, we can see that
\[ \{\phi_{\alpha^{k-1}(x_0)}\dots\phi_{\alpha(x_0)}\phi_{x_0}(t):
k=m_{n-1},2m_{n-1},\dots,l_nm_{n-1}\} \]
is $1/n$-dense in $\T$.
Hence, for every $t\in\T$, the closure of
\[ \{(\alpha\times\phi)^k(x_0,t):k\in\N\} \]
in $X\times\T$ contains $\{x_0\}\times\T$.
It follows that $\alpha\times\phi$ is minimal.
Let us check that $\phi$ is never cohomologous to
a cocycle with values in the rotation group.
Suppose that $\phi$ is cohomologous to a cocycle with values
in the rotation group.
Then we would have a homeomorphism $\id\times\gamma$ on $X\times\T$
such that
\[ (\id\times\gamma)(\alpha\times\phi)=
(\alpha\times\phi)(\id\times\gamma) \]
and $(\id\times\gamma)(x_0,0)=(x_0,1/2)$.
By (3), we can verify that
\[ (\alpha\times\phi)^{m_n}(x_0,0)\rightarrow(x_0,0) \]
and
\[ (\alpha\times\phi)^{m_n}(x_0,1/2)\rightarrow(x_0,0) \]
in $X\times\T$ as $n\rightarrow\infty$.
Consequently we obtain
\begin{align*}
(\id\times\gamma)(x_0,0)
&=\lim_{n\rightarrow\infty}
(\id\times\gamma)(\alpha\times\phi)^{m_n}(x_0,0) \\
&=\lim_{n\rightarrow\infty}
(\alpha\times\phi)^{m_n}(\id\times\gamma)(x_0,0) \\
&=\lim_{n\rightarrow\infty}
(\alpha\times\phi)^{m_n}(x_0,1/2)=(x_0,0),
\end{align*}
which is a contradiction.

Let us construct $m_n$, $l_n$ and $\phi_n$.
Put $m_1=l_1=1$ and let $\phi_1$ be an irrational rotation.
Suppose that $m_{n-1}$, $k_{n-1}$ and $\phi_{n-1}$ have been fixed.
Put
\[ \psi_{n-1}=\phi_{n-1}(m_{n-1}-1)\dots\phi_{n-1}(1)\phi_{n-1}(0). \]
By (2), there exists $\omega\in\Homeo^+(\T)$ such that
$\omega\psi_{n-1}\omega^{-1}=R_\theta$,
where $r(\psi_{n-1})=\theta$ is an irrational number.
Since $\psi_{n-1}$ is minimal,
we can find a natural number $l_n$ so that the condition (4) holds.
Applying Fact (a) to $\omega(0)$ and $\omega(1/2)$,
we obtain $\rho\in\Homeo^+(\T)$ such that
both $\lvert\omega^{-1}\rho\omega(0)\rvert$ and
$\lvert\omega^{-1}\rho\omega(1/2)\rvert$ are less than $1/n$.
Furthermore $\rho$ is conjugate to an irrational rotation.
Define $\tilde{\omega}:\Z_{m_{n-1}}\rightarrow\Homeo^+(\T)$ by
\[ \tilde{\omega}(k)=
\omega\phi_{n-1}(0)^{-1}\phi_{n-1}(1)^{-1}\dots\phi_{n-1}(k-1)^{-1} \]
for all $k=0,1,\dots,m_{n-1}-1$.
Evidently we have
\[ \tilde{\omega}(k+1)\phi_{n-1}(k)\tilde{\omega}(k)^{-1}=
\begin{cases} \id & k\neq m_{n-1}-1 \\
R_\theta & k=m_{n-1}-1 \end{cases} \]
for all $k\in\Z_{m_{n-1}}$.
Choose $\varepsilon>0$ so that $d(\phi,\id)<\varepsilon$ implies
\[ d(\phi_{n-1}(k)\tilde{\omega}(k)^{-1}\phi\tilde{\omega}(k),
\phi_{n-1}(k))<\frac{1}{2^n} \]
for all $k\in\Z_{m_{n-1}}$.
By Fact (b), we can find a natural number $N$ greater than $l_n$
and a sequence of homeomorphisms
\[ \id=\rho_0,\rho_1,\dots,\rho_{m_{n-1}(N-l_n)}=
\rho R_\theta^{-N} \]
such that
\[ d(\rho_{i+1}\rho_i^{-1},\id)<\varepsilon \]
for all $i=0,1,\dots,m_{n-1}(N-l_n)-1$.
Note that this is easily done because $d(\cdot,\cdot)$ is
invariant under rotations.
Put $m_n=m_{n-1}N$.
We define a map $\tilde{\rho}:\Z_{m_n}\rightarrow\Homeo^+(\T)$ by
\[ \tilde{\rho}(k)=
\begin{cases}
\id & k=0,1,\dots,l_nm_{n-1}-1 \\
R_\theta^{-j}\rho_{k'+1}\rho_{k'}^{-1}R_\theta^{j} & \text{otherwise},
\end{cases}\]
where $k'=k-l_nm_{n-1}$ and $j$ is a natural number satisfying
$N-k/m_{n-1}\leq j<N+1-k/m_{n-1}$.
Notice that we still have
\[ d(\tilde{\rho}(k),\id)<\varepsilon \]
for all $k\in\Z_{m_n}$.
Define $\phi_n:\Z_{m_n}\rightarrow\Homeo^+(\T)$ by
\[ \phi_n(k)=\phi_{n-1}(k)
\tilde{\omega}(k)^{-1}\tilde{\rho}(k)\tilde{\omega}(k). \]
The condition (5) is already built in this definition.
The condition (6) is immediate from the choice of $\varepsilon$.
Since one can check that
\begin{align*}
\psi_n&=\phi_n(m_n-1)\dots\phi_n(1)\phi_n(0) \\
&=\tilde{\omega}(m_n)^{-1}\left(
R_\theta\tilde{\rho}(m_n-1)\dots\tilde{\rho}(m_n-m_{n-1})
R_\theta\dots\tilde{\rho}(l_nm_{n-1})
\right)R_\theta^{l_n}\tilde{\omega}(0) \\
&=\omega^{-1}\rho R_\theta^{-N}R_\theta^N\omega \\
&=\omega^{-1}\rho\omega,
\end{align*}
the conditions (2) and (3) follow immediately.
\end{exm}

\flushleft{
\textit{Huaxin Lin\\
e-mail: hxlin@noether.uoregon.edu\\
Department of Mathematics\\
University of Oregon\\
Eugene, Oregon 97403\\
U.S.A.\\}}

\flushleft{
\textit{Hiroki Matui\\
e-mail: matui@math.s.chiba-u.ac.jp \\
Graduate School of Science and Technology,\\
Chiba University,\\
1-33 Yayoi-cho, Inage-ku,\\
Chiba 263-8522,\\
Japan. }}

\end{document}